\documentclass[twocolumn]{autart}    

\usepackage{cite}
\usepackage{amsmath,amssymb,amsfonts}
\usepackage{algorithmic}
\usepackage{algorithm}
\usepackage[bookmarks=false]{hyperref}
\usepackage{textcomp}
\usepackage{subfigure}
\usepackage{float}
\usepackage{placeins}

\usepackage{graphicx}
\newtheorem{theorem}{Theorem}[section]

\newtheorem{lemma}[theorem]{Lemma}
\newtheorem{proposition}[theorem]{Proposition}
\newtheorem{corollary}[theorem]{Corollary}
\newtheorem{assumption}[theorem]{Assumption}
\newtheorem{remark}[theorem]{Remark}

\newenvironment{proof}{%
  \par\noindent\textbf{Proof.}\ }{%
  \hfill$\square$\par}

\newcommand{\rev}[1]{#1}
\newcommand{\bluerev}[1]{#1}

\begin{document}

\begin{frontmatter}

\title{\bluerev{Linear--Quadratic Mean Field Games under Heterogeneous
Erroneous Initial Information}}

\thanks[footnoteinfo]{This work was supported by the National Natural Science Foundation of China (NSFC Grant No.12441101).}
\thanks{Acknowledgements: The authors thank Prof. Markus Fischer and Prof. Alekos Cecchin for helpful discussions.}
\author[SYHC,SMS,LMIB]{Yuxin Jin}\ead{yxjin@buaa.edu.cn}
\author[CAICT]{Lu Ren}\ead{urrenlu@buaa.edu.cn}
\author[SAI,LMIB]{Wang Yao\corauthref{cor1}}\ead{yaowang@buaa.edu.cn}
\author[SMS,LMIB]{Xiao Zhang\corauthref{cor1}}\ead{xiao.zh@buaa.edu.cn}
\corauth[cor1]{Corresponding authors: Wang Yao and Xiao Zhang.}
\address[SYHC]{ShenYuan Honors College, Beihang University, Beijing 100191, China}
\address[SMS]{School of Mathematical Sciences, Beihang University, Beijing 100191, China}
\address[LMIB]{Key Laboratory of Mathematics, Informatics and Behavioral Semantics, Ministry of Education, Beihang University, Beijing 100191, China}
\address[SAI]{School of Artificial Intelligence, Beihang University, Beijing 100191, China}
\address[CAICT]{China Academy of Information and Communications Technology, Beijing 100191, China}


\begin{abstract}
\bluerev{We study a finite-horizon linear--quadratic mean field game with
heterogeneous observations of the initial mean field. These observations may
be erroneous and are used by agents to construct their feedback laws. At the
population level, the heterogeneous error profile is generally
infinite-dimensional. We derive exact linear sensitivity formulas showing
that its closed-loop effects nevertheless admit a finite-dimensional closure:
for each agent, the deviations are governed by only two \(n\)-dimensional
error channels, its private error \(E_i\) and the population-average error
\(\bar E\). The population-average error determines the displacement of the
actual mean field, whereas centered private errors determine agent-specific
deviations and the additional cross-sectional covariance.} The same
representation yields quadratic cost identities and, under the stated
i.i.d.\ finite-population construction, a uniform \(O(N^{-1})\) mean-square
approximation of the empirical aggregate. In the deterministic model, the
private state history induces a finite-dimensional linear observation problem
for the private error and the population-average error. Nonsingularity of the
associated observability Gramian is necessary and sufficient for exact
recovery. The recovered errors reconstruct the actual mean-field state at the
revision time and initialize a revised correct-information continuation from
that state. We obtain an exact continuation-cost comparison; when the
population-average error vanishes, revision does not increase the continuation
cost. For stochastic dynamics, we take a two-component revision signal as
given and study the resulting one-shot feedback relative to an oracle
continuation initialized with the actual current mean field. The actual
mean-field deviation is driven only by the population-average signal error, while
individual deviations also retain the private signal error. Exact covariance
and quadratic continuation-cost identities quantify these effects.
\end{abstract}
\end{frontmatter}

\section{Introduction}

Mean field games (MFGs) provide a tractable framework for large-population dynamic games by representing population interactions through mean-field quantities in a limiting model \cite{1,2,3,4,5,6}. In linear--quadratic (LQ) settings, the equilibrium analysis often reduces to Riccati equations and linear consistency systems \cite{7,8,9,10,genlq}. A distinctive feature of this setting is that the mean field plays two roles simultaneously: it enters the feedback calculation of each agent and is itself generated by the population of resulting feedbacks. An error in the initial mean-field observation can therefore propagate through the closed-loop interaction, acting more than a perturbation of an individual initial state.

This paper studies this effect when the initial mean-field observations are heterogeneous. An agent \(\mathcal A_i\) observes $z_0^i=z_0+E_i$ and treats this value as the correct initial mean field state $z_0$, computes its subjective mean field prediction and selects its feedback at the initial time $t=0$. Such feedbacks generate the actual mean-field evolution. This raises three questions. First, how do heterogeneous initial information errors propagate through the private and population dynamics? Second, can an agent infer the relevant errors from its own state trajectory? Third, if agents can change feedbacks at a revision time \(t_0\), how will the game evolve?

Related works consider several information-constrained MFG settings. Partially observed MFGs derive controls from filtered or restricted state information \cite{poCommon,21,22}, while quantized-information models study the effect of limited observation resolution \cite{20}. Adaptive MFGs allow unknown model quantities to be learned during the evolution \cite{23}, and incomplete-information or Bayesian-learning models introduce richer belief structures \cite{bertucci2022,shmayaZiliotto}. Related LQG control problems study estimation and control under communication constraints \cite{comm}. A complementary implementation study \cite{confcomp} develops a
method based on discrete private-state observations and repeated control
updates. 

\bluerev{The present paper focuses on the structural propagation and
identifiability of heterogeneous initial mean-field observation errors. At the
population level, these errors form a generally infinite-dimensional profile,
yet their linear propagation through the affine LQ system admits an exact
finite-dimensional closure. For any fixed agent \(\mathcal A_i\), the
closed-loop deviations are driven by only two \(n\)-dimensional error channels:
its private error \(E_i\) and the population-average error \(\bar E\). The
latter determines the displacement of the actual mean field, whereas
centered private errors determine deviations around that aggregate. This
closure also turns deterministic trajectory-based recovery into a
finite-dimensional linear observability problem
\cite{kalman1960,kailath1980}.}

The main contributions are as follows. 
\begin{itemize} 
\item We derive exact linear representations for the propagation of heterogeneous initial mean-field observation errors. The population-average error determines the displacement of the actual mean field, while centered private errors determine deviations around the actual mean field and the corresponding covariance correction. The same representation yields quadratic-cost identities. For a finite population formed from i.i.d.\ representative copies, we further obtain a uniform \(O(N^{-1})\) mean-square approximation of the empirical aggregate by its limiting actual mean field. 
\item In the deterministic model, the private state history over \([0,t_0]\) encodes both the private error and the population-average error through a finite-dimensional linear observation family. We show that nonsingularity of the associated observability Gramian is necessary and sufficient for exact recovery of \((E_i^\top,\bar E^\top)^\top\) on the unrestricted error space. The recovered quantities reconstruct the actual mean-field state \(z_{t_0}^A\). The revised continuation is then initialized from this current state, rather than from the original benchmark trajectory, so that revision corrects the information while preserving the actual mean-field state generated before \(t_0\). We derive the resulting state, covariance, and continuation-cost relations. For a general biased population, the cost difference is parameter dependent; when \(\bar E=0\), revision does not increase the continuation cost. 
\item We then consider stochastic dynamics with an exogenous two-component signal available at the revision time. The signal specifies an estimate of the current mean field and an estimate of the population-average estimate. It selects an affine feedback that is optimal for the corresponding deterministic subjective environment and is held fixed on \([t_0,T]\). Relative to an oracle continuation initialized with the actual current mean field, the resulting actual mean-field deviation depends only on the population-average signal error, whereas the tagged-agent deviation also retains the private signal error. Exact covariance and quadratic continuation-cost identities quantify the resulting first- and second-order effects. 
\end{itemize} 
The correct-information LQ-MFG equilibrium is the benchmark. The numerical experiments illustrate the propagation
formulas, finite-population scaling, deterministic revision, and
continuation-cost effects for the chosen parameters.

The rest of this paper is organized as follows: Section~2 introduces the model and the correct-information LQ-MFG benchmark.
Section~3 derives the error-propagation, covariance, quadratic-cost, and
finite-population relations. Section~4 develops deterministic identification
and exact revision. Section~5 studies signal-based revision under a
two-component signal, and Section~6 provides numerical illustrations of the
theoretical results.
\section{Problem formulation and correct-information benchmark}\label{sec:model}
We begin with the finite-population LQ model and then collect the
correct-information benchmark objects used in Sections~3--5.

\subsection{Finite-population model and information structure}
\rev{Let $(\Omega,\mathcal F,\mathbb F,\mathbb P)$ be a complete filtered
probability space satisfying the usual conditions, and fix $T>0$.}
We take
\[
\begin{gathered}
A,C,\Gamma,\bar\Gamma\in\mathbb R^{n\times n},
B,F\in\mathbb R^{n\times d},D\in\mathbb R^{n\times r_W},\\
s,\bar s,\eta,\bar\eta\in\mathbb R^n,
W^i\text{ is }\mathbb R^{r_W}\text{-valued},\\
Q_I,Q,\bar Q_I,\bar Q\in\mathbb R^{n\times n},
R\in\mathbb R^{d\times d}.
\end{gathered}
\]
For each $1\le i\le N$, $\mathcal A_i$ has state process $x^i\in\mathbb R^n$ governed by
\begin{equation}
\left\{
\begin{aligned}
&dx^i_t=\big(Ax^i_t+Bu^i_t+Cx^{(N)}_t+F\bar u^{(N)}_t\big)dt+DdW^i_t,\\
&x^i(0)=x_0^i,
\end{aligned}
\right.
\label{eq:sec2_dyn_N}
\end{equation}
where $u^i$ is the control of $\mathcal{A}_i$,
\[
x^{(N)}_t:=\frac1N\sum_{j=1}^N x^j_t,
\bar u^{(N)}_t:=\frac1N\sum_{j=1}^N u^j_t,
\]
and \(\{W^i\}_{i=1}^N\) are independent standard
$\mathbb F$-Brownian motions in $\mathbb R^{r_W}$, while
$\{x_0^i\}_{i=1}^N$ are $\mathcal F_0$-measurable square-integrable initial
states.
The cost functional of $\mathcal A_i$ is
\begin{equation}
\label{eq:sec2_cost_N}
\left\{
\begin{aligned}
&J_i(u^i,u^{-i})=\frac12\mathbb E\Bigl[\int_0^T\bigl(\|x^i_t-s\|_{Q_I}^2+\|u^i_t\|_{R}^2\\
&+\|x^i_t-(\Gamma x^{(N)}_t+\eta)\|_{Q}^2\bigr)dt+\|x^i_T-\bar s\|_{\bar Q_I}^2\\
&+\|x^i_T-(\bar\Gamma x^{(N)}_T+\bar\eta)\|_{\bar Q}^2\Bigr].
\end{aligned}
\right.
\end{equation}
Here \(\|y\|_M^2:=y^\top M y\).
All model matrices are deterministic and constant, \(T<\infty\), and the
generic initial state satisfies \(X_0\in L^2(\Omega;\mathbb R^n)\). We assume
\[
Q_I,Q,\bar Q_I,\bar Q\in\mathbb S_+^n,
R\in\mathbb S_{++}^d.
\]
Set \(K_B:=R^{-1}B^\top\).

Let \(\Theta\) denote the deterministic model parameters listed above. The
initial mean \(z_0\) and the realized error variables are not included in
\(\Theta\).

Individual optimization uses state-feedback controls. For an
interval \(I=[a,b]\subseteq[0,T]\), deterministic paths \((z,\bar u)\), and a
prescribed square-integrable initial condition, let
\[
\mathfrak U^{\mathrm{fb}}_I(z,\bar u)
\]
denote the Borel maps
\(\alpha:I\times\mathbb R^n\to\mathbb R^d\) for which the closed-loop equation
\[
dx_t=\bigl(Ax_t+B\alpha(t,x_t)+Cz_t+F\bar u_t\bigr)dt+DdW_t
\]
has a unique strong solution satisfying
\[
\mathbb E\sup_{t\in I}|x_t|^2
+\mathbb E\int_a^b|\alpha(t,x_t)|^2\mathrm dt<\infty.
\]
The feedback is selected using the information available at the beginning of
\(I\) and is held fixed on that interval. The same well-posedness and
integrability requirements apply to finite-population profiles in
\eqref{eq:sec2_dyn_N}. All individual optimality statements below concern the
limiting problem with deterministic mean-field paths.

\rev{Under correct initial information, the limiting control problem is
parameterized by $\Theta$ and the initial mean field}
\[
z_0:=\mathbb E[X_0],
\]
\rev{where \(X_0\) is a generic population state. }

\subsection{Correct-information benchmark and feedback equilibrium}
\rev{In the mean-field limit, the correct-information benchmark uses
deterministic paths \(z,\bar u\) and the affine adjoint
\(p_t^i=P_t^1x_t^i+g_t\) for the LQ-MFG defined by}
\eqref{eq:sec2_dyn_N}--\eqref{eq:sec2_cost_N}; see
\cite{1,2,3,7,8,9,10}. The individual Riccati matrix \(P^1\) solves
\begin{equation}
\left\{
\begin{aligned}
&-\dot P^1_t=P^1_tA+A^\top P^1_t+(Q_I+Q)-P^1_tBK_BP^1_t,\\
&P^1_T=\bar Q_I+\bar Q,
\end{aligned}
\right.
\label{eq:sec2_P1}
\end{equation}
and $g$ solves
\begin{equation}
\left\{
\begin{aligned}
&\dot g_t=-((A^\top-P^1_tBK_B)g_t+(P^1_tC-Q\Gamma)z_t\\
&+P^1_tF\bar u_t-Q_Is-Q\eta),\\
&g_T=-\bar Q_I\bar s-\bar Q(\bar\Gamma z_T+\bar\eta).
\end{aligned}
\right.
\label{eq:sec2_g}
\end{equation}
where \(z\) and \(\bar u\) are the prescribed mean-field state and control
paths.
The benchmark feedback is
\begin{equation}
u^i_t=\varphi^c(t,x^i_t):=-K_B\big(P^1_tx^i_t+g_t\big).
\label{eq:sec2_phi_c}
\end{equation}
Set
\[
\mathcal Q:=Q_I+Q-Q\Gamma,\bar{\mathcal Q}:=\bar Q_I+\bar Q-\bar Q\bar\Gamma.
\]
Then the deterministic mean field forward--backward system is
\begin{equation}
\left\{
\begin{aligned}
&\dot z_t=(A+C)z_t-(B+F)K_Bp_t,\ z(0)=z_0,\\
&\dot p_t=-(A^\top p_t+\mathcal Q z_t-Q_I s-Q\eta),\\
&p_T=\bar{\mathcal Q}z_T-\bar Q_I\bar s-\bar Q\bar\eta.
\end{aligned}
\right.
\label{eq:sec2_mf_fbsde}
\end{equation}
Equivalently, writing
\[p_t=P^0_tz_t+\mathcal G_t,\]
we obtain
\begin{equation}
\dot z_t
=
(A+C-(B+F)K_BP^0_t)z_t
-(B+F)K_B\mathcal G_t,
\label{eq:sec2_z}
\end{equation}
where $P^0$ solves
\begin{equation}
\left\{
\begin{aligned}
&-\dot P^0_t=P^0_t(A+C)+A^\top P^0_t-P^0_t(B+F)K_BP^0_t\\
&+\mathcal Q,P^0_T=\bar{\mathcal Q},
\end{aligned}
\right.
\label{eq:sec2_P0}
\end{equation}
and \(\mathcal G\) solves
\begin{equation}
\left\{
\begin{aligned}
&\dot{\mathcal G}_t=-[(A^\top-P^0_t(B+F)K_B)\mathcal G_t-Q_I s-Q\eta],\\
&\mathcal G_T=-\bar Q_I\bar s-\bar Q\bar\eta.
\end{aligned}
\right.
\label{eq:sec2_Gaffine}
\end{equation}
\rev{Take the solution of \eqref{eq:sec2_z} with initial value \(z_0\) as
\(z^c\), and set}
\begin{equation}
\bar u_t^c:=-K_B(P_t^0z_t^c+\mathcal G_t),t\in[0,T].
\label{eq:sec2_ubar_c}
\end{equation}
\rev{The pair \((z^c,\varphi^c)\) defined by
\eqref{eq:sec2_z}--\eqref{eq:sec2_Gaffine} and
\eqref{eq:sec2_phi_c} is the \emph{correct-information benchmark}. It is
individually optimal and mean-field consistent within the stated deterministic
mean-field and affine state-feedback class.
The affine-adjoint verification is given in the Supplementary Material.}

We use \(\bar u^c\), \(g^c\), and \(p^c\) for the corresponding benchmark
mean control, individual adjoint offset, and mean-field adjoint.

\begin{remark}\label{rem:g_affine_z}

For any self-consistent affine continuation on \([\tau,T]\),
\(0\leq\tau<T\), whose averaged state and adjoint satisfy the benchmark
mean-field forward--backward system, the adjoint term in the individual
feedback satisfies
\[
g_t=(P_t^0-P_t^1)z_t+\mathcal G_t .
\]
\end{remark}
\rev{Averaging the individual adjoint equation gives
\(g_t=(P_t^0-P_t^1)z_t+\mathcal G_t\); details are provided in the
Supplementary Material.}

\subsection{Benchmark solvability and population convention}
\emph{Standing LQ conditions.}
Under \(Q_I,Q,\bar Q_I,\bar Q\in\mathbb S_+^n\) and
\(R\in\mathbb S_{++}^d\), the individual Riccati
equation \eqref{eq:sec2_P1} has its standard unique bounded solution
\(P^1\in C^1([0,T];\mathbb S^n)\). We additionally assume that the generally
nonsymmetric mean-field Riccati equation \eqref{eq:sec2_P0} admits a unique
bounded solution \(P^0\in C^1([0,T];\mathbb R^{n\times n})\) on \([0,T]\).
\rev{Related finite-horizon LQ mean-field solvability analyses using Riccati
or forward--backward systems can be found in \cite{7,dom,delay,genlq}.}

\rev{Given \(P^0\), \eqref{eq:sec2_Gaffine} and
\eqref{eq:sec2_z} have unique solutions for every prescribed initial mean.}

\section{Error propagation in LQMFGs under heterogeneous erroneous initial
information}\label{sec:errorprop}

In this section, we answer the question: how initial errors affects the system? Starting from the correct-information benchmark \((z^c,\varphi^c)\), we derive linear maps from the initial errors to the deviations in subjective mean-field
predictions, feedback laws, actual mean field, individual trajectories,
population covariance, and cost.

\begin{assumption}
\label{ass:sec3_probabilistic}
The generic initial state $X_0$ has finite second moment, and $W$
is a standard Brownian motion independent of $X_0$. The generic
initial observation error $E$ is independent of $X_0$ and $W$, and has finite second moment.
\end{assumption}
Recall \(z_0=\mathbb E[X_0]\), and define its mean and covariance by
\[
\bar E:=\mathbb E[E],
\Sigma^E:=\mathbb E[(E-\bar E)(E-\bar E)^\top].
\]

\subsection{Decentralized response under erroneous initial mean field observations}

\begin{assumption}
\label{ass:sec3_erroneous_observations}
At time \(t=0\), \(\mathcal A_i\) receives the deterministic value
\[
z_0^i = z_0 + E_i,
\]
where \(E_i\in\mathbb R^n\) is its realized initial mean-field observation
error. 
\end{assumption}

\begin{assumption}
\label{ass:sec3_subjective_response}
\(\mathcal A_i\) uses \(z_0^i\) as the initial mean field and keeps the
resulting subjective continuation fixed on \([0,T]\).
\end{assumption}

\rev{Let \(z^i\) and \(\bar u^i\) denote the deterministic subjective
mean-field state and control. Given these paths, \(\mathcal A_i\) solves}
\[
\inf_{\alpha\in\mathfrak U^{\mathrm{fb}}_{[0,T]}(z^i,\bar u^i)}
J_i(\alpha;z^i,\bar u^i),
\]
\rev{where \(J_i(\alpha;z^i,\bar u^i)\) is the functional in
\eqref{eq:sec2_cost_N} with \(u_t^i=\alpha(t,x_t^i)\),
\(x_t^{(N)}\) replaced by \(z_t^i\), and \(\bar u_t^{(N)}\) replaced by
\(\bar u_t^i\). The feedback is computed from
\((\Theta,z_0^i)\) at the initial time.}

Under an admissible feedback \(\alpha\), with
\(u_t^i=\alpha(t,x_t^i)\), the subject private state of $\mathcal A_i$ satisfies
\[
\left\{
\begin{aligned}
&dx_t^i=\bigl(Ax_t^i+Bu_t^i+Cz_t^i+F\bar u_t^i\bigr)dt+DdW_t^i,\\
&x^i(0)=x_0^i.
\end{aligned}
\right.
\]
The LQ structure gives the feedback representation
\begin{equation}
u_t^i = \varphi^i(t,x_t^i) := -K_B\bigl(P^1_{t}x_t^i + g_t^i\bigr),
\label{eq:sec3_feedback_err}
\end{equation}
The matrices \(P^0_t,P^1_t\) and the affine term \(\mathcal G_t\) are the
benchmark coefficients determined by the common model parameters. The
subjective mean field prediction is determined by
\begin{equation}
\left\{
\begin{aligned}
&\dot z_t^i=(A+C-(B+F)K_BP_{t}^0)z_t^i-(B+F)K_B\mathcal G_t,\\
&z_0^i=z_0+E_i,
\end{aligned}
\right.
\label{eq:sec3_zi}
\end{equation}
and the subjective mean field control is
\begin{equation}
\bar u_t^i = -K_B\bigl(P^0_{t} z_t^i + \mathcal G_t\bigr).
\label{eq:sec3_ubari}
\end{equation}
The auxiliary term $g^i$ satisfies
\begin{equation}
\left\{
\begin{aligned}
&\dot g_t^i=-[(A^\top-P^1_{t}BK_B)g_t^i+(P^1_{t}C-Q\Gamma)z_t^i+P^1_{t}F\bar u_t^i\\
&-Q_Is-Q\eta],g_T^i=-\bar Q_I\bar s-\bar Q(\bar\Gamma z_T^i+\bar\eta).
\end{aligned}
\right.
\label{eq:sec3_gi}
\end{equation}

For a perceived initial mean field \(y\), let
\(\mathfrak S_0(\Theta,y)\) denote the feedback obtained from
\eqref{eq:sec3_zi}--\eqref{eq:sec3_gi} with \(z_0^i=y\). Thus
\[
\varphi^i=\mathfrak S_0(\Theta,z_0^i).
\]

\subsection{Linear propagation of initial errors}
\rev{Compare the subjective objects \((z^i,\bar u^i,g^i)\) and the actual
evolution with the Section~2 benchmark, and define}
\[
\delta z_t^i := z_t^i - z_t^c,
\delta \bar u_t^i := \bar u_t^i - \bar u_t^c,
\delta g_t^i := g_t^i - g_t^c.
\]
Let \(x_t^{i,c}\) and \(x_t^{i,A}\) use the same initial state \(x_0^i\) and
Brownian path \(W^i\).

For the generic error \(E\), let
\((z^E,g^E,\varphi^E)\) denote the subjective objects generated from
\(z_0+E\), equivalently
\(\varphi^E:=\mathfrak S_0(\Theta,z_0+E)\). A generic agent with mark
\((X_0,E,W)\) uses the control \(u_t=\varphi^E(t,X_t^A)\). Substituting this
feedback into the limiting individual dynamics yields the actual state
\begin{equation}\label{eq:sec3_background_dynamics}
\left\{
\begin{aligned}
&dX_t^A=\bigl[AX_t^A+B\varphi^E(t,X_t^A)+Cz_t^A+F\bar u_t^A\bigr]dt+\\
&DdW_t, X_0^A=X_0,z_t^A=\mathbb E[X_t^A],\bar u_t^A=\mathbb E[\varphi^E(t,X_t^A)].
\end{aligned}
\right.
\end{equation}
The identities \(z_t^A=\mathbb E[X_t^A]\) and
\(\bar u_t^A=\mathbb E[\varphi^E(t,X_t^A)]\) impose mean-field consistency for
the profile \(\varphi^E\).
Write \(\delta g_t^E:=g_t^E-g_t^c\).
Let \(X^c\) be the benchmark state driven by the same initial condition
\(X_0\) and Brownian motion \(W\), with feedback \(\varphi^c\) and aggregates
\((z^c,\bar u^c)\).
The dynamics of $\mathcal A_i$ under the actual mean field and mean control
\((z^A,\bar u^A)\) are
\[
\left\{
\begin{aligned}
&dx_t^{i,A}=\bigl[Ax_t^{i,A}+B\varphi^i(t,x_t^{i,A})+Cz_t^A+F\bar u_t^A\bigr]dt\\
&+DdW_t^i,x_0^{i,A}=x_0^i .
\end{aligned}
\right.
\]
Also set
\[
u_t^{i,A}:=\varphi^i(t,x_t^{i,A}),
u_t^{i,c}:=\varphi^c(t,x_t^{i,c}).
\]
Define the deviations of the actual mean field, tagged-agent state, and mean
control by
\[
\begin{aligned}
\delta z_t^A:=z_t^A-z_t^c,\delta x_t^i:=x_t^{i,A}-x_t^{i,c},\delta\bar u_t^A:=\bar u_t^A-\bar u_t^c.
\end{aligned}
\]
All deviations below are linear maps of the private error \(E_i\) and the
population mean error \(\bar E\). The fundamental matrices encode homogeneous
closed-loop propagation, while \(M^g,M^z,M^{x,1},M^{x,2}\) record how these
two error channels enter the feedback offset, actual mean field, and
individual state.

Let \(\Phi^1,\Phi^z,\Phi^x\) be the fundamental matrices
\begin{equation}
\dot \Phi_t^1
=(A+C-(B+F)K_BP_{t}^0)\Phi_t^1,
\Phi_0^1 = I_n.
\label{eq:sec3_phi1}
\end{equation}
\begin{equation}
\dot \Phi_t^z=
(A+C-(B+F)K_BP_{t}^1)\Phi_t^z,
\Phi_0^z = I_n.
\label{eq:sec3_phiz}
\end{equation}
\begin{equation}
\dot \Phi_t^x
=(A-BK_BP_{t}^1)\Phi_t^x,
\Phi_0^x = I_n.
\label{eq:sec3_phix}
\end{equation}
Define
\begin{equation}\label{eq:sec3_Mg}
M_t^g:=(P_t^0-P_t^1)\Phi_t^1 .
\end{equation}
\begin{equation}\label{eq:sec3_Mz}
M_t^z:=-\Phi_t^z\int_0^t(\Phi_s^z)^{-1}(B+F)K_BM_s^g\mathrm{d}s.
\end{equation}
\begin{equation}\label{eq:sec3_Mx}
\begin{aligned}
M_t^{x,1}&=-\Phi_t^x\int_0^t(\Phi_s^x)^{-1}BK_BM_s^g\mathrm{d}s,\\
M_t^{x,2}&=M_t^z-M_t^{x,1}.
\end{aligned}
\end{equation}
Set
\begin{equation}\label{eq:sec3_error_to_deviation_maps}
\left\{
\begin{aligned}
&\theta_i:=\begin{pmatrix}E_i\\\bar E\end{pmatrix},X_t:=\begin{pmatrix}M_t^{x,1}&M_t^{x,2}\end{pmatrix},Z_t:=\begin{pmatrix}0&M_t^z\end{pmatrix},\\
&G_t^E:=\begin{pmatrix}M_t^g&0\end{pmatrix},U_t:=-K_B(P_t^1X_t+G_t^E).
\end{aligned}
\right.
\end{equation}

\begin{theorem}\label{thm:sec3_linear_propagation}

Under the standing LQ conditions and
Assumptions~\ref{ass:sec3_probabilistic},
\ref{ass:sec3_erroneous_observations}, and
\ref{ass:sec3_subjective_response}, for \(t\in[0,T]\),
\begin{equation}\label{eq:sec3_deltazi_representation}
\delta z_t^i=\Phi_t^1E_i,
\end{equation}
\begin{equation}\label{eq:sec3_deltaubari_representation}
\delta\bar u_t^i=-K_BP_t^0\Phi_t^1E_i,
\end{equation}
\begin{equation}\label{eq:sec3_deltagi_representation}
\delta g_t^i=M_t^gE_i,
\end{equation}
\begin{equation}\label{eq:sec3_deltazA_representation}
\delta z_t^A=M_t^z\bar E,
\end{equation}
\begin{equation}\label{eq:sec3_deltaubarA_representation}
\delta\bar u_t^A=-K_B(P_t^1M_t^z+M_t^g)\bar E,
\end{equation}
and
\begin{equation}\label{eq:sec3_deltax_sol}
\delta x_t^i=M_t^{x,1}E_i+M_t^{x,2}\bar E.
\end{equation}
Equivalently, with \(\delta u_t^i:=u_t^{i,A}-u_t^{i,c}\),
\[
\begin{aligned}
\delta x_t^i=X_t\theta_i,
\delta z_t^A=Z_t\theta_i,\delta g_t^i=G_t^E\theta_i,
\delta u_t^i=U_t\theta_i.
\end{aligned}
\]
\end{theorem}
\begin{proof}
Subtracting the benchmark equations from the subjective equations first gives
the representations of \(\delta z^i\), \(\delta\bar u^i\), and
\(\delta g^i\). Averaging the generic closed-loop deviation then yields a
linear equation for \(\delta z^A\) driven only by \(\bar E\). Variation of constants
produces the stated matrices, and the control representation follows from the
affine feedback law. The complete proof is given in the Supplementary
Material.
\end{proof}
\rev{The actual mean field and mean control depend on the error law only
through \(\bar E\), the aggregate bias channel. The centered channel
\(E_i-\bar E\) drives an agent's deviation from the actual mean field, while
\(E-\bar E\) generates cross-sectional dispersion.}

For the covariance calculation, let
\[
\Sigma_0:=\mathbb E[(X_0-z_0)(X_0-z_0)^\top],
\]
and define
\[
\left\{
\begin{aligned}
&\Sigma_t^c:=\mathbb E[(X_t^c-z_t^c)(X_t^c-z_t^c)^\top],\\
&\Sigma_t^A:=\mathbb E[(X_t^A-z_t^A)(X_t^A-z_t^A)^\top].
\end{aligned}
\right.
\]

\begin{proposition}\label{prop:sec3_covariance}
Under the assumptions of Theorem~\ref{thm:sec3_linear_propagation}, for every
\(t\in[0,T]\),
\begin{equation}\label{eq:sec3_prop_sigmac}
\Sigma_t^c=
\Phi_t^x\left(\Sigma_0+\int_0^t(\Phi_s^x)^{-1}DD^\top(\Phi_s^x)^{-\top}\mathrm{d}s\right)(\Phi_t^x)^\top.
\end{equation}
The covariance under erroneous observations is
\begin{equation}\label{eq:sec3_covariance_identity}
\Sigma_t^A=\Sigma_t^c+M_t^{x,1}\Sigma^E(M_t^{x,1})^\top .
\end{equation}
\end{proposition}
\rev{Thus \(\Sigma_t^A\) is the cross-sectional covariance of the limiting
population law: the mean error shifts the actual mean field, whereas centered
private errors add the positive semidefinite correction in
\eqref{eq:sec3_covariance_identity}.}

\subsection{Finite-population implementation}

\rev{For a finite-\(N\) implementation, let \(x^{j,N,A}\) denote the
closed-loop state from \eqref{eq:sec2_dyn_N} when coordinate \(j\) uses
\(\varphi^{E_j^N}:=\mathfrak S_0(\Theta,z_0+E_j^N)\):}
\[
\left\{
\begin{aligned}
&dx_t^{j,N,A}=\bigl[Ax_t^{j,N,A}
+B\varphi^{E_j^N}(t,x_t^{j,N,A})+Cx_t^{(N),A}\\
&+F\bar u_t^{(N),A}\bigr]dt
+DdW_t^{j,N},x_0^{j,N,A}=x_0^{j,N},\\
&x_t^{(N),A}:=\frac1N\sum_{k=1}^Nx_t^{k,N,A},\\
&\bar u_t^{(N),A}:=\frac1N\sum_{k=1}^N
\varphi^{E_k^N}(t,x_t^{k,N,A}).
\end{aligned}
\right.
\]
Set
\[
\left\{
\begin{aligned}
&x_0^{(N)}:=N^{-1}\sum_{j=1}^Nx_0^{j,N},
E^{(N)}:=N^{-1}\sum_{j=1}^NE_j^N,\\
&W_t^{(N)}:=N^{-1}\sum_{j=1}^NW_t^{j,N}.
\end{aligned}
\right.
\]
\begin{proposition}\label{prop:sec3_finite_population}
For the finite-population profile defined above, in which agent \(j\) uses
\(\varphi^{E_j^N}\),
\begin{equation}\label{eq:sec3_prop37}
\begin{aligned}
x_t^{(N),A}-z_t^A
=&\Phi_t^z(x_0^{(N)}-z_0)+M_t^z(E^{(N)}-\bar E)\\
&+\Phi_t^z\int_0^t(\Phi_s^z)^{-1}DdW_s^{(N)}.
\end{aligned}
\end{equation}
\rev{If \((x_0^{j,N},E_j^N,W^{j,N})_{j=1}^N\) are i.i.d. copies of the generic
background triple in Assumption~\ref{ass:sec3_probabilistic}, then a constant
\(C_T>0\) independent of \(N\) satisfies}
\begin{equation}\label{eq:sec3_finite_ms}
\sup_{0\le t\le T}\mathbb E|x_t^{(N),A}-z_t^A|^2\le\frac{C_T}{N}.
\end{equation}
\end{proposition}

\subsection{Cost distortion under erroneous initial observations}
For \(\mathcal A_i\), fix \(E_i\) and \(\bar E\). If \(x_0^i\) is fixed, the
expectation below is over the Brownian noise; if \(x_0^i\) is random, it is
included as well. The control deviation satisfies
\begin{equation}\label{eq:sec3_cost_ui}
\delta u_t^i=-K_B(P_t^1\delta x_t^i+\delta g_t^i).
\end{equation}
Set
\begin{equation}\label{eq:sec3_JA}
\left\{
\begin{aligned}
&J_i^A:=\frac12\mathbb E\Bigl[\int_0^T(\|x_t^{i,A}-s\|_{Q_I}^2+\|u_t^{i,A}\|_R^2\\
&+\|x_t^{i,A}-(\Gamma z_t^A+\eta)\|_Q^2)\mathrm{d}t+\|x_T^{i,A}-\bar s\|_{\bar Q_I}^2\\
&+\|x_T^{i,A}-(\bar\Gamma z_T^A+\bar\eta)\|_{\bar Q}^2\Bigr],\\
&J_i^c:=\frac12\mathbb E\Bigl[\int_0^T(\|x_t^{i,c}-s\|_{Q_I}^2+\|u_t^{i,c}\|_R^2\\
&+\|x_t^{i,c}-(\Gamma z_t^c+\eta)\|_Q^2)\mathrm{d}t+\|x_T^{i,c}-\bar s\|_{\bar Q_I}^2\\
&+\|x_T^{i,c}-(\bar\Gamma z_T^c+\bar\eta)\|_{\bar Q}^2\Bigr].
\end{aligned}
\right.
\end{equation}

The linear deviation identities imply a quadratic cost representation in
\(\theta_i\). The Supplementary Material defines its coefficients and proves
the centered specialization. The continuation costs used in Section~4 are
defined below.

\section{Error identification and strategy revision in the deterministic case}

\rev{Throughout this section, we answer the question: when \(D=0\), how can an agent identify initial errors, recover the actual system state and revise its strategy? Fix
\(t_0\in(0,T)\). Conditional on the initial state and realized errors, the
trajectories are deterministic. We first derive the linear observation system.
When its Gramian is nonsingular, the private history up to \(t_0\) identifies
\((E_i,\bar E)\), reconstructs \(z_{t_0}^A\), and initializes the revised
continuation, denoted by the superscript \(r\).}

\rev{The deviation identities in Section~3 reduce recovery to a
finite-dimensional linear inverse problem.}
\subsection{Information from the private state trajectory}
Set
\[
\Psi_x(t,s):=\Phi_t^x(\Phi_s^x)^{-1},0\le s\le t\le t_0.
\]
Variation of constants separates the terms determined by the subjective path
\((z^i,g^i)\) and known model coefficients from the two unknown error
components.
For \(t\in[0,t_0]\), define the integral observation
\begin{equation}\label{eq:sec4_integral_observation}
\begin{aligned}
Y_t^i:=&x_t^{i,A}-\Psi_x(t,0)x_0^i-\int_0^t\Psi_x(t,s)\{(C-FK_BP_s^1)z_s^i\\
&-(B+F)K_Bg_s^i\}\mathrm{d}s .
\end{aligned}
\end{equation}
\rev{The quantity \(Y_t^i\) is computed from the private state trajectory, the
initial state, and the subjective paths \(z^i,g^i\) generated by \(z_0^i\).}

Define
\[
\left\{
\begin{aligned}
K_t^E&:=(FK_BP_t^0-C)\Phi_t^1,\\
K_t^{\bar E}&:=(C-FK_BP_t^1)M_t^z-FK_BM_t^g,\\
K_t&:=\begin{pmatrix}K_t^E&K_t^{\bar E}\end{pmatrix},
\xi_i:=\begin{pmatrix}E_i\\\bar E\end{pmatrix}.
\end{aligned}
\right.
\]
Variation of constants and Theorem~\ref{thm:sec3_linear_propagation} give
\begin{equation}\label{eq:sec4_integral_observation_linear}
\begin{aligned}
Y_t^i&=H_t\xi_i,H_t:=\int_0^t\Psi_x(t,s)K_s\mathrm{d}s\in\mathbb R^{n\times2n}.
\end{aligned}
\end{equation}
\subsection{Identifiability of the initial errors}
The family \(t\mapsto H_t\xi_i\) is a deterministic linear observation system
for \(\xi_i\). Its observability over \([0,t_0]\) is measured by the integral
Gramian
\begin{equation}\label{eq:sec4_WH}
\mathcal W_H(t):=\int_0^tH_s^\top H_s\mathrm{d}s,0\le t\le t_0.
\end{equation}

\begin{theorem}\label{thm:sec4_continuous_recoverability}
On the parameter space \(\mathbb R^{2n}\), the
full vector \(\xi_i=(E_i^\top,\bar E^\top)^\top\) is uniquely recoverable from
\(\{Y_t^i:0\le t\le t_0\}\) if and only if \(\mathcal W_H(t_0)\) is
nonsingular. In that case,
\begin{equation}\label{eq:sec4_integral_reconstruction}
\begin{aligned}
\widehat\xi_i:=\mathcal W_H(t_0)^{-1}
\int_0^{t_0}H_s^\top Y_s^i\mathrm{d}s
=\xi_i,\xi_i=\begin{pmatrix}E_i\\\bar E\end{pmatrix}.
\end{aligned}
\end{equation}
Moreover, if
\[
t_*:=\inf\{t\in(0,t_0]:\det\mathcal W_H(t)\neq0\}
\]
exists, then \(\mathcal W_H(t)\) is nonsingular for every
\(t\in(t_*,t_0]\), and
\[
\mathcal D_t(x^i_{[0,t]},z_0^i)
:=\mathcal W_H(t)^{-1}\int_0^tH_s^\top Y_s^i\,\mathrm ds
=\xi_i .
\]
\end{theorem}
\begin{proof}
If \(\mathcal W_H(t_0)\) is nonsingular and two candidates
\(\xi,\tilde\xi\) generate the same observations, then
\(H_t(\xi-\tilde\xi)=0\) for all \(t\in[0,t_0]\). Hence
\[
(\xi-\tilde\xi)^\top\mathcal W_H(t_0)(\xi-\tilde\xi)
=\int_0^{t_0}|H_s(\xi-\tilde\xi)|^2\mathrm{d}s=0,
\]
so \(\xi=\tilde\xi\). The reconstruction formula follows from
\[
\int_0^{t_0}H_s^\top Y_s^i\mathrm{d}s
=\int_0^{t_0}H_s^\top H_s\mathrm{d}s\xi_i
=\mathcal W_H(t_0)\xi_i .
\]
Conversely, if \(\mathcal W_H(t_0)\) is singular, then some nonzero
\(\delta\) satisfies
\[
\int_0^{t_0}|H_s\delta|^2\mathrm{d}s=0 .
\]
By continuity, \(H_t\delta=0\) on \([0,t_0]\), so \(\xi_i\) and
\(\xi_i+\delta\) generate the same integral observations. Unique recovery is
therefore impossible. The final assertion follows from the monotonicity
\(\mathcal W_H(t)=\mathcal W_H(s)+\int_s^tH_\tau^\top H_\tau d\tau\) for
\(t\ge s\).
\end{proof}

\begin{remark}
The matrix \(\mathcal W_H(t_0)\) is the observability Gramian of the linear
observation family \(Y_t^i=H_t\xi_i\). Its nonsingularity means that the
time family \(H_t\) separates the private error and population-average error
channels.
Through mean-field coupling, the private trajectory encodes both channels;
the Gramian tests whether their time signatures are distinguishable.
\end{remark}

\rev{Because \(H_t\) is determined by \(\Theta\), the Gramian can be evaluated
before revision. The subsequent analysis assumes exact recoverability and uses
\(\widehat\xi_i=\mathcal D_{t_0}(x^i_{[0,t_0]},z_0^i)\).}

\begin{assumption}\label{ass:exact-recoverability}
The revision time \(t_0\) is chosen such that the continuous-time Gramian
\(\mathcal W_H(t_0)\) is nonsingular.
\end{assumption}

\subsection{Reconstruction of the actual mean field state at the revision time}
\begin{proposition}\label{prop:sec4_reconstruct_zA}
\rev{Under Assumption~\ref{ass:exact-recoverability}, $\mathcal A_i$ reconstructs
the actual mean-field state at time $t_0$ through}
\begin{equation}
z_{t_0}^A=z_{t_0}^i+M_{t_0}^z\bar E-\Phi_{t_0}^1E_i.
\end{equation}
\rev{In particular, the reconstructed value is common to all agents.}
\end{proposition}

\begin{proof}
\rev{Section~3 gives}
\[
z_t^A-z_t^i=M_t^z\bar E-\Phi_t^1E_i .
\]
Setting \(t=t_0\) gives the displayed reconstruction formula, whose value is
the common actual mean-field state \(z_{t_0}^A\).
\end{proof}

\subsection{Revised feedback law}
\rev{At the common time \(t_0\), every population member applies the decoder
to its private history, reconstructs the same \(z_{t_0}^A\), and switches
simultaneously to the revised correct-information continuation initialized at
that state.}

Let $z_t^r$ and $\bar u_t^r$ be defined by
\begin{equation}
\left\{
\begin{aligned}
&\dot z_t^r=(A+C-(B+F)K_BP_t^0)z_t^r-(B+F)K_B\mathcal G_t,\\
&z_{t_0}^r=z_{t_0}^A,
\end{aligned}
\right.
\end{equation}
and
\begin{equation}
\bar u_t^r=-K_B(P_t^0z_t^r+\mathcal G_t),t\in[t_0,T].
\end{equation}
Define $g_t^r$ on $[t_0,T]$ by
\begin{equation}
\left\{
\begin{aligned}
&\dot g_t^r=-\Big[(A^\top-P_t^1BK_B)g_t^r
+(P_t^1C-Q\Gamma)z_t^r\\
&+P_t^1F\bar u_t^r-Q_Is-Q\eta\Big],\\
&g_T^r=-\bar Q_I\bar s-\bar Q(\bar\Gamma z_T^r+\bar\eta).
\end{aligned}
\right.
\end{equation}
The revised feedback law is
\begin{equation}
u_t^{i,r}=\varphi^r(t,x_t^{i,r}):=-K_B(P_t^1x_t^{i,r}+g_t^r),t\in[t_0,T],
\end{equation}
and the corresponding revised state trajectory satisfies
\begin{equation}
\left\{
\begin{aligned}
&\dot x_t^{i,r}=A x_t^{i,r}+B u_t^{i,r}+C z_t^r+F \bar u_t^r,\\
&x_{t_0}^{i,r}=x_{t_0}^{i,A}.
\end{aligned}
\right.
\end{equation}
Let \(X^r\) denote the generic background counterpart of \(x^{i,r}\), with
initial condition \(X_{t_0}^A\). Its population mean and mean control are
\(z^r\) and \(\bar u^r\).

\rev{By Remark~\ref{rem:g_affine_z},
\(g_t^r=(P_t^0-P_t^1)z_t^r+\mathcal G_t\). Averaging the generic closed-loop
equation under \(\varphi^r\) yields the mean state \(z^r\) and mean control
\(\bar u^r\); the verification is given in the Supplementary Material.}

\subsection{Deviations after revision}
Set
\[
\begin{aligned}
&\delta z_t^r:=z_t^r-z_t^c,\delta x_t^{i,r}:=x_t^{i,r}-x_t^{i,c}.
\end{aligned}
\]
\rev{For \(t,s\in[t_0,T]\), define \(\Psi_1\) and use
\(\Psi_x(t,s)=\Phi_t^x(\Phi_s^x)^{-1}\):}
\[
\Psi_1(t,s):=\Phi_t^1(\Phi_s^1)^{-1}.
\]

\begin{theorem}
\label{prop:sec4_post_revision_deviation}
For \(t\in[t_0,T]\),
\begin{equation}\label{eq:sec4_post_revision_state}
\left\{
\begin{aligned}
&\delta z_t^r=\Psi_1(t,t_0)M_{t_0}^z\bar E,\\
&\delta x_t^{i,r}=\Psi_x(t,t_0)M_{t_0}^{x,1}E_i+\bigl[\Psi_1(t,t_0)M_{t_0}^z\\
&-\Psi_x(t,t_0)M_{t_0}^{x,1}\bigr]\bar E.
\end{aligned}
\right.
\end{equation}

Moreover, for
\(\Sigma_t^r:=\operatorname{Cov}(X_t^r-z_t^r)\), 
\begin{equation}\label{eq:sec4_post_revision_covariance}
\Sigma_t^r
=
\Psi_x(t,t_0)\Sigma_{t_0}^A\Psi_x(t,t_0)^\top .
\end{equation}
\end{theorem}

\begin{proof}
The equations for \(z^r\) and \(z^c\) have the same homogeneous part, and
\(\delta z_{t_0}^r=M_{t_0}^z\bar E\). This gives the first identity in
\eqref{eq:sec4_post_revision_state}. Also
\[
g_t^r-g_t^c=(P_t^0-P_t^1)\delta z_t^r .
\]
Therefore
\[
\begin{aligned}
&\delta\dot x_t^{i,r}
=(A-BK_BP_t^1)\delta x_t^{i,r}+A_t^{xz}\delta z_t^r,\\
&A_t^{xz}:=C-(B+F)K_BP_t^0+BK_BP_t^1 .
\end{aligned}
\]
Using variation of constants and
\[
\frac{d}{ds}\bigl[\Psi_x(t,s)\Psi_1(s,t_0)\bigr]
=\Psi_x(t,s)A_s^{xz}\Psi_1(s,t_0),
\]
we obtain
\[
\delta x_t^{i,r}
=\Psi_x(t,t_0)\delta x_{t_0}^{i,A}
+[\Psi_1(t,t_0)-\Psi_x(t,t_0)]\delta z_{t_0}^A .
\]
Substituting
\[
\delta x_{t_0}^{i,A}
=M_{t_0}^z\bar E+M_{t_0}^{x,1}(E_i-\bar E),
\delta z_{t_0}^A=M_{t_0}^z\bar E
\]
gives \eqref{eq:sec4_post_revision_state}. Subtracting the first line of
\eqref{eq:sec4_post_revision_state} from the second gives
\[
X_t^r-z_t^r=\Psi_x(t,t_0)(X_{t_0}^A-z_{t_0}^A),
\]
and the covariance identity follows.
\end{proof}

\subsection{Continuation cost comparison}

For \(k\in\{A,r,c\}\), let \(J_{i,t_0}^k\) be the objective continuation cost
on \([t_0,T]\) evaluated along \((x^{i,k},u^{i,k},z^k)\):
\begin{equation}\label{eq:sec4_continuation_cost}
\left\{
\begin{aligned}
&J_{i,t_0}^k:=\frac12\mathbb E\Bigl[\int_{t_0}^T\bigl(\|x_t^{i,k}-s\|_{Q_I}^2+\|u_t^{i,k}\|_R^2\\
&+\|x_t^{i,k}-(\Gamma z_t^k+\eta)\|_Q^2\bigr)\mathrm{d}t+\|x_T^{i,k}-\bar s\|_{\bar Q_I}^2\\
&+\|x_T^{i,k}-(\bar\Gamma z_T^k+\bar\eta)\|_{\bar Q}^2\Bigr].
\end{aligned}
\right.
\end{equation}
\rev{For \(k\in\{A,r\}\), let
\(\ell_{i,t_0}^k\) and \(\mathcal H_{t_0}^k\) denote the linear and quadratic
coefficients obtained by substituting the corresponding linear deviations into
\eqref{eq:sec4_continuation_cost}. Their integral definitions are given in the
Supplementary Material.}

\begin{theorem}\label{thm:sec4_performance_gap}
Under the standing LQ conditions and
Assumption~\ref{ass:exact-recoverability},
for \(\mathcal A_i\), fixed \(E_i\), and fixed \(\bar E\),
\begin{equation}\label{eq:sec4_exact_performance_gap}
J_{i,t_0}^A-J_{i,t_0}^r
=
(\ell_{i,t_0}^A-\ell_{i,t_0}^r)^\top\theta_i
+\frac12\theta_i^\top(\mathcal H_{t_0}^A-\mathcal H_{t_0}^r)\theta_i .
\end{equation}
\end{theorem}
\rev{The identity follows by substituting the linear state, mean-field, and
control deviations into the quadratic cost; the details are given in the
Supplementary Material.}

\rev{For a general biased population, the sign of the right-hand side of
\eqref{eq:sec4_exact_performance_gap} depends on the model parameters and the realized errors.} When $\bar{E}=0$, the right-hand side is nonnegative.
\begin{corollary}\label{cor:sec4_centered_performance_gain}
Suppose that \(\bar E=0\). Then
\(z^A=z^r=z^c\) and \(\bar u^A=\bar u^r=\bar u^c\) on \([t_0,T]\), and
\[
J_{i,t_0}^A-J_{i,t_0}^r
=
\frac12\int_{t_0}^T
E_i^\top(M_t^g)^\top BK_BM_t^gE_i\mathrm{d}t\ge0.
\]
Equality holds if and only if
\[
\begin{aligned}
B^\top M_t^gE_i=0,\text{for a.e. }t\in[t_0,T].
\end{aligned}
\]
\end{corollary}
\begin{proof}
When $\bar E=0$, the erroneous and revised continuations evolve in the same
deterministic mean-field environment, and the revised feedback is optimal for
that environment. Applying the finite-horizon completion-of-squares identity on
\([t_0,T]\), with initial state \(x_{t_0}^{i,A}\), gives
\[
J_{i,t_0}^A-J_{i,t_0}^r
=\frac12\mathbb E\int_{t_0}^T
\bigl\|u_t^{i,A}-\varphi^r(t,x_t^{i,A})\bigr\|_R^2\mathrm dt.
\]
Since $g_t^i-g_t^r=M_t^gE_i$ and both feedbacks have the same state
coefficient,
\[
u_t^{i,A}-\varphi^r(t,x_t^{i,A})=-K_BM_t^gE_i.
\]
Using $K_B=R^{-1}B^\top$ and the symmetry of $R$ yields the displayed
quadratic integral. Because $R\succ0$, equality holds exactly when
$B^\top M_t^gE_i=0$ for almost every $t$.
\end{proof}

\rev{When \(\bar E=0\), both controls face the same mean-field environment;
revision removes the feedback-offset error induced by the private error
\(g_t^i-g_t^r=M_t^gE_i\), which yields the definite sign above.}
\setlength{\emergencystretch}{2em}

\section{Signal-based revision under a two-component signal}
\label{sec:two_layer_revision}
\rev{Section~4 treats deterministic recovery and revision. Here a given
two-component signal at \(t_0\) selects a feedback law that remains fixed on
\([t_0,T]\).
We define the augmented-information and restricted two-layer constructions,
then derive exact signal-error propagation relative to an oracle continuation.}

\[
\begin{aligned}
&M_B:=(B+F)K_B,\\
&L^k_t:=A+C-M_BP^k_t,k\in\{0,1\}.
\end{aligned}
\]

\subsection{Two-component signal and signal-based feedback}
\begin{assumption}
\label{ass:sec5_signal_device}
At \(t_0\), \(\mathcal A_i\) observes its current state and the two-component
signal
\[
\sigma_i=(\bar\zeta_i,\hat\zeta_i)\in\mathbb R^n\times\mathbb R^n .
\]
\rev{The component \(\hat\zeta_i\) is \(\mathcal A_i\)'s estimate of the current
actual mean-field state \(z_{t_0}^A\). The component \(\bar\zeta_i\) is its
estimate of the population average of the current mean-field estimates.}
\bluerev{\(W_t-W_{t_0}\) and \(W_t^i-W_{t_0}^i\), \(t\in[t_0,T]\), remain Brownian
motions.}
\end{assumption}

\begin{assumption}
\label{ass:sec5_one_shot_admissibility}
After observing \(\sigma_i\), \(\mathcal A_i\) selects a Borel feedback
map
\[
\psi_{\sigma_i}:[t_0,T]\times\mathbb R^n\to\mathbb R^d, \tilde u_t=\psi_{\sigma_i}(t,x_t).
\]
The selected map belongs to
\bluerev{\(\mathfrak U^{\mathrm{fb}}_{[t_0,T]}
(z^{\sigma_i},\bar u^{\sigma_i})\)}, where
\bluerev{\((z^{\sigma_i},\bar u^{\sigma_i})\)} are defined by
\eqref{eq:sec5_restricted_average_environment} and
\eqref{eq:sec5_hat_environment}. The map is selected at \(t_0\) and thereafter
evaluated only at the current state.
\end{assumption}

\rev{For aggregation, let \(\mathsf S=(\bar\zeta,\hat\zeta)\) have law
\(\nu_{\mathsf S}:=\mathcal L(\mathsf S)\), jointly with
\((E,X^A_{[0,t_0]})\), and set}
\[
m_{\hat\zeta}:=\mathbb E[\hat\zeta].
\]
\rev{Here \(m_{\hat\zeta}\) is the population mean of the current estimates;
\(\sigma_i\), \(\sigma\), and \(\mathsf S\) denote the realized, fixed, and
generic signals, respectively.}

\begin{assumption}
\label{ass:sec5_subjective_common}
\bluerev{When computing \(\psi_{\sigma_i}\), \(\mathcal A_i\) treats
\(\hat\zeta_i\) as the current mean-field state and \(\bar\zeta_i\) as the
population mean of the current estimates. In its own subjective calculation,
it therefore substitutes \(m_{\hat\zeta}=\bar\zeta_i\). This is a belief used
to close the agent's calculation; it does not assert objectively that
\(\bar\zeta_i=\bar\zeta_j\) across agents.}
\end{assumption}
\bluerev{This assumption closes the agent's subjective system.}

\subsection{Two-layer closure}
\bluerev{Consider an agent who receives a fixed signal
\(\sigma=(\bar\zeta,\hat\zeta)\). Write \(z^\sigma,\bar u^\sigma\) for its
predicted actual mean field state and actual mean field control,
\(x^\sigma,u^\sigma\) for its state and control, and \(g^\sigma\) for the
adjoint term in its feedback control. Write \(\bar z,\bar g\) as the
population average prediction of the mean field state and the population
average adjoint term in the feedback control, and \(z^\nu,\bar u^\nu\) as the
actual signal-based mean field and mean control. We first consider an augmented
information case, where agents share \(\bar z\) and \(\bar g\).}

\noindent\textit{Augmented information.}
\bluerev{Assume that agents share \(\bar g\) and \(\bar z\). For each fixed
signal value \(\sigma\) and corresponding \(z^\sigma\), the subjective optimal
response is}
\[
u_t^{\sigma}=-K_B(P_t^1x_t^{\sigma}+g_t^{\sigma}),
\]
\bluerev{where}
\begin{equation}
\left\{
\begin{aligned}
&\bluerev{dx_t^\sigma}
 =\big(Ax_t^\sigma+Bu_t^\sigma+C\bluerev{z_t^\nu}
 +F\bluerev{\bar u_t^\nu}\big)\bluerev{dt}
 +\bluerev{D\,dW_t^\sigma},x_{t_0}^\sigma=x_{t_0}^\sigma,\\
&\dot g_t^\sigma=-(A^\top-P_t^1BK_B)g_t^\sigma
 -(P_t^1C-Q\Gamma)\bluerev{z_t^\sigma}\\
&-P_t^1F\bar u^\sigma_t+Q_Is+Q\eta,\\
&g_T^\sigma=-\bar Q_I\bar s
 -\bar Q(\bar\Gamma\bluerev{z_T^\sigma}+\bar\eta).
\end{aligned}
\right.
\label{eq:sec5_augmented_response}
\end{equation}
\bluerev{The first line keeps the current tautological initial condition only
to display the response in the standard state-adjoint form. Since
\(z_t^\nu=\mathbb E[x_t^{\mathsf S}]\),
\(\bar g_t=\mathbb E[g_t^{\mathsf S}]\), and
\(\bar u_t^\nu=\mathbb E[u_t^{\mathsf S}]\), combining with
Assumption~\ref{ass:sec5_subjective_common}, we have}
\[
\bluerev{\bar u_t^\nu}=-K_B(P_t^1\bluerev{z_t^\nu}+\bar g_t),
\bar u_t^\sigma=-K_B(P_t^1z_t^\sigma+\bar g_t),
\]
\bluerev{and}
\begin{equation}\label{eq:sec5_aug_zg}
\left\{
\begin{split}
&\dot{\bluerev{z_t^\nu}}
=(A+C-(B+F)K_BP_t^1)\bluerev{z_t^\nu}-(B+F)K_B\bar g_t,\\
&\dot z_t^\sigma
=(A+C-(B+F)K_BP_t^1)z_t^\sigma-(B+F)K_B\bar g_t,\\
&\dot{\bar z}_t
=(A+C-(B+F)K_BP_t^1)\bar z_t-(B+F)K_B\bar g_t,\\
&z_{t_0}^\sigma=\hat\zeta,\bar z_{t_0}=m_{\hat\zeta},\\
&\dot{\bar g}_t
=-(A^\top-P_t^1(B+F)K_B)\bar g_t-(P_t^1C-Q\Gamma\\
&-P_t^1FK_BP_t^1)\bar z_t+Q_Is+Q\eta,\\
&\bar g_T=-\bar Q_I\bar s-\bar Q(\bar\Gamma\bar z_T+\bar\eta).
\end{split}
\right.
\end{equation}
\bluerev{Thus, for a given \(\sigma\), \(\bar z,\bar g,z^\sigma\) can be
computed by solving \eqref{eq:sec5_aug_zg}, and the feedback can be evaluated.}

\bluerev{Then we consider the restricted information case, where an agent only
knows its own signal.}

\noindent\textit{Restricted information.}
\bluerev{Under restricted information, according to
Assumption~\ref{ass:sec5_subjective_common}, \(\mathcal A_i\) uses
\(m_{\hat\zeta}=\bar\zeta_i\) in its own calculation. Hence \(\mathcal A_i\)
can solve the average layer \(\bar z^{\sigma_i}\) from \(\bar\zeta_i\):}
\begin{equation}\label{eq:sec5_restricted_average_environment}
\left\{
\begin{split}
&\dot{\bar z}^{\sigma_i}_t
=(A+C-(B+F)K_BP_t^1)\bar z^{\sigma_i}_t
-(B+F)K_B\bar g^{\sigma_i}_t,\\
&\bar z^{\sigma_i}_{t_0}=\bar\zeta_i,\\
&\dot{\bar g}^{\sigma_i}_t
=-(A^\top-P_t^1(B+F)K_B)\bar g^{\sigma_i}_t-(P_t^1C-Q\Gamma\\
&-P_t^1FK_BP_t^1)\bar z^{\sigma_i}_t+Q_Is+Q\eta,\\
&\bar g^{\sigma_i}_T=-\bar Q_I\bar s
-\bar Q(\bar\Gamma\bar z^{\sigma_i}_T+\bar\eta).
\end{split}
\right.
\end{equation}
\bluerev{Then \(\mathcal A_i\) defines
\(\bar u_t^{\sigma_i}:=-K_B(P_t^1z_t^{\sigma_i}+\bar g_t^{\sigma_i})\) and
solves the response layer \(z^{\sigma_i},g^{\sigma_i}\) from
\(\hat\zeta_i\):}
\begin{equation}\label{eq:sec5_hat_environment}
\left\{
\begin{split}
&\dot z_t^{\sigma_i}
=(A+C-(B+F)K_BP_t^1)z_t^{\sigma_i}
-(B+F)K_B\bar g^{\sigma_i}_t,\\
&z_{t_0}^{\sigma_i}=\hat\zeta_i,\\
&\dot g_t^{\sigma_i}=-(A^\top-P_t^1BK_B)g_t^{\sigma_i}
-(P_t^1C-Q\Gamma)\bluerev{z_t^{\sigma_i}}\\
&-P_t^1F\bar u^{\sigma_i}_t+Q_Is+Q\eta,\\
&g_T^{\sigma_i}=-\bar Q_I\bar s
-\bar Q(\bar\Gamma\bluerev{z_T^{\sigma_i}}+\bar\eta).
\end{split}
\right.
\end{equation}
\bluerev{The resulting subjective optimal feedback is}
\begin{equation}\label{eq:sec5_restricted_two_layer}
\psi_{\sigma_i}(t,x)=-K_B(P_t^1x+\bluerev{g_t^{\sigma_i}}).
\end{equation}
\bluerev{Once \((z^{\sigma_i},\bar u^{\sigma_i})\) is frozen, this is the same
deterministic-environment LQ problem and the same optimality argument as in
Sections~2--3. }
\subsection{Propagation of signal errors}
\label{subsec:sec5_no_update_profile}

\rev{Given \(\sigma_i=(\bar\zeta_i,\hat\zeta_i)\), agent \(\mathcal A_i\)
solves \eqref{eq:sec5_restricted_average_environment}--
\eqref{eq:sec5_restricted_two_layer} and
\bluerev{keeps \(\psi_{\sigma_i}\) fixed on \([t_0,T]\)}. Set}
\[
\bluerev{g_t^{i,\nu}:=g_t^{\sigma_i},\qquad
\varphi^{i,\nu}(t,x):=\psi_{\sigma_i}(t,x).}
\]
\rev{The superscript \(\nu\) denotes the actual continuation in which each
population member uses the feedback associated with its own signal. Its
population mean state and mean control are}
\[
z_t^{\nu}:=\mathbb E[X_t^{\nu}],
\bluerev{\bar u_t^{\nu}:=\mathbb E[\psi_{\mathsf S}(t,X_t^{\nu})]},
\]
where
\[
\left\{
\begin{aligned}
&dX_t^{\nu}
=\big(A X_t^{\nu}+B\bluerev{\psi_{\mathsf S}(t,X_t^{\nu})}\big)dt+\big(Cz_t^{\nu}+F\bar u_t^{\nu}\big)dt\\
&+DdW_t,X_{t_0}^{\nu}=X_{t_0}^A .
\end{aligned}
\right.
\]

Let \(z^\star:=z_{t_0}^A\). In the stochastic oracle continuation, the
mean-field state is initialized at \(z^\star\), while each individual state
retains its actual value at the revision time. Define the continuation by
\[
\left\{
\begin{aligned}
&\dot z_t^o=(A+C-(B+F)K_BP_t^0)z_t^o-(B+F)K_B\mathcal G_t,\\
&z_{t_0}^o=z^\star,\\
&\bar u_t^o=-K_B(P_t^0z_t^o+\mathcal G_t),g_t^o=(P_t^0-P_t^1)z_t^o+\mathcal G_t,\\
&\varphi^o(t,x)=-K_B(P_t^1x+g_t^o).
\end{aligned}
\right.
\]
The generic and \(\mathcal A_i\)'s oracle states satisfy
\[
\left\{
\begin{aligned}
&dX_t^o=\bigl(AX_t^o+B\varphi^o(t,X_t^o)\bigr)dt+\bigl(Cz_t^o+F\bar u_t^o\bigr)dt]\\
&+DdW_t,X_{t_0}^o=X_{t_0}^A,\\
&dx_t^{i,o}=\bigl(Ax_t^{i,o}+B\varphi^o(t,x_t^{i,o})\bigr)dt+\bigl(Cz_t^o+F\bar u_t^o\bigr)dt\\
&+DdW_t^i,x_{t_0}^{i,o}=x_{t_0}^{i,A},\\
&u_t^{i,o}:=\varphi^o(t,x_t^{i,o}).
\end{aligned}
\right.
\]
The oracle consistency identities are
\[
\mathbb E[X_t^o]=z_t^o,
\mathbb E[\varphi^o(t,X_t^o)]=\bar u_t^o.
\]
The oracle and revised continuations use the same Brownian increments after
\(t_0\). If \(D=0\), exact recovery gives
\[
(z^o,\bar u^o,g^o,\varphi^o,x^{i,o})
=(z^r,\bar u^r,g^r,\varphi^r,x^{i,r}).
\]

The state and control of \(\mathcal A_i\) under this signal-based profile satisfy
\[
\left\{
\begin{aligned}
&dx_t^{i,\nu}
=\big(Ax_t^{i,\nu}+Bu_t^{i,\nu}+Cz_t^{\nu}+F\bar u_t^{\nu}\big)dt+DdW_t^i,\\
&u_t^{i,\nu}=\varphi^{i,\nu}(t,x_t^{i,\nu}),
x_{t_0}^{i,\nu}=x_{t_0}^{i,o}=x_{t_0}^{i,A},
\end{aligned}
\right.
\]

\rev{First-moment identities require \(e^{\mathsf S}\in L^1\), whereas
covariance and quadratic-cost identities require
\(e^{\mathsf S}\in L^2\).}

For the realized signal
\(\sigma_i=(\bar\zeta_i,\hat\zeta_i)\), set
\[
e_i^{\bar z}:=\bar\zeta_i-z^\star,
e_i^{\hat z}:=\hat\zeta_i-z^\star,
\]
and define
\[
\begin{aligned}
&e_i:=\begin{pmatrix}e_i^{\bar z}\\e_i^{\hat z}\end{pmatrix},
e^{\mathsf S}:=\begin{pmatrix}\bar\zeta-z^\star\\ \hat\zeta-z^\star\end{pmatrix},\bar e:=\mathbb E[e^{\mathsf S}].
\end{aligned}
\]

Here \(e_i^{\hat z}\) is the current mean-field error. The first component is
measured relative to \(z^\star\), whereas its deviation from the population
target \(m_{\hat\zeta}\) is \(\bar\zeta_i-m_{\hat\zeta}\).

{
\begin{theorem}
\label{thm:sec5_signal_error_propagation}
Under the standing LQ conditions and
Assumptions~\ref{ass:sec5_signal_device},
\ref{ass:sec5_one_shot_admissibility}, and
\ref{ass:sec5_subjective_common}, assume
\(e^{\mathsf S}\in L^1\). There is a deterministic continuous sensitivity
matrix \(\mathsf C_t\in\mathbb R^{n\times 2n}\), depending only on
\((\Theta,t_0)\), such that
\[
g_t^{i,\nu}-g_t^o=\mathsf C_te_i.
\]
Moreover,
\begin{equation}
z_t^\nu-z_t^o
=-\Phi_t^z\int_{t_0}^t
(\Phi_s^z)^{-1}M_B\mathsf C_s\bar e\,ds,
\label{eq:sec5_mean_error_direct}
\end{equation}
and, under the same-noise coupling,
\begin{equation}
\begin{aligned}
x_t^{i,\nu}-x_t^{i,o}=&\Phi_t^x\int_{t_0}^t(\Phi_s^x)^{-1}\big[(C-FK_BP_s^1)(z_s^\nu-z_s^o)\\
&-FK_B\mathsf C_s\bar e-BK_B\mathsf C_se_i\big]ds.
\end{aligned}
\label{eq:sec5_state_error_direct}
\end{equation}
The control difference is therefore
\begin{equation}
u_t^{i,\nu}-u_t^{i,o}
=-K_B\bigl[P_t^1(x_t^{i,\nu}-x_t^{i,o})
+\mathsf C_te_i\bigr].
\label{eq:sec5_control_error_direct}
\end{equation}
In particular, the actual mean-field deviation depends only on the
population-average signal error \(\bar e\), whereas the individual state also
contains the private signal error \(e_i\).
\end{theorem}
The construction of \(\mathsf C\), together with the centered-state covariance
and direct continuation-cost identities, is given in the Supplementary
Material.
}
\subsection{Diagonal single-layer special case}
\label{subsec:sec5_single_layer_special}

\bluerev{If the two estimates coincide,}
\[
\bar\zeta_i=\hat\zeta_i=\zeta_i,
\]
\bluerev{forward uniqueness first gives equality of the two state predictions;
with this common path, backward uniqueness gives equality of the two adjoint
terms and hence of the predicted mean controls:}
\[
\bluerev{\bar z^{\sigma_i}=z^{\sigma_i}=:z^{e,i},\qquad
\bar g^{\sigma_i}=g^{\sigma_i}=:g^{e,i}.}
\]
\bluerev{The two-layer construction therefore collapses to}
\begin{equation}
\left\{
\begin{aligned}
&\dot z^{e,i}_t=L^1_tz^{e,i}_t-M_Bg^{e,i}_t,z^{e,i}_{t_0}=\zeta_i,\\
&\dot g^{e,i}_t
=-\Big[(A^\top-P^1_tM_B)g^{e,i}_t+(P^1_tC-P^1_tFK_BP^1_t\\
&-Q\Gamma)z^{e,i}_t-Q_Is-Q\eta\Big],\\
&g^{e,i}_T=-\bar Q_I\bar s-\bar Q(\bar\Gamma z^{e,i}_T+\bar\eta),\\
&\varphi^{e,i}(t,x)=-K_B(P^1_tx+g^{e,i}_t),u^{e,i}_t=\varphi^{e,i}(t,x_t).
\end{aligned}
\right.
\label{eq:sec5_single_layer_p1}
\end{equation}
\bluerev{The affine identity
\(g^{e,i}_t=(P_t^0-P_t^1)z^{e,i}_t+\mathcal G_t\) also gives uniqueness of
this single-layer continuation. Section~6 uses only this diagonal
specialization; the general two-layer construction is retained for the
signal-error propagation and covariance formulas.}


\section{Numerical experiments}\label{sec:numerics}
\rev{We use deterministic and stochastic experiments to illustrate propagation
of the actual mean field, finite-population approximation, exact revision, and
signal-based revision.}

\subsection{Experimental setup}
We use the following parameter set
\[
\begin{aligned}
&T=2, \Delta t=T/400,t_0=0.5,A=-I_2, B=0.5I_2, \\
&C=I_2, F=0.5I_2,R=I_2, Q_I=Q=\bar Q_I=\bar Q=I_2,\\
&\Gamma=\bar\Gamma=I_2, \eta=\bar\eta=0,s=\bar s=(0.5,0.3)^\top,\\
&z_0=(0.5,-0.5)^\top,\bar E=(0.5,-0.5)^\top,\\
&\Sigma^E=0.04I_2,\operatorname{Cov}(X_0)=0.03I_2 .
\end{aligned}
\]
Initial states and observation errors are sampled independently from Gaussian
laws with the displayed moments and independently of the Brownian motion.
For trajectory plots, sampled initial states and signal errors are centered at
their prescribed means. The signal-based revision uses
\[
D=0.05I_2,
\bar\varepsilon=(0.1,-0.1)^\top,
\Sigma^\varepsilon=0.01I_2,
\]
and the diagonal single-layer signal in
Subsection~\ref{subsec:numerics_single_layer}.

\subsection{Linear propagation of the actual mean field}
\begin{figure}[t]
    \centering
    \includegraphics[width=\linewidth]{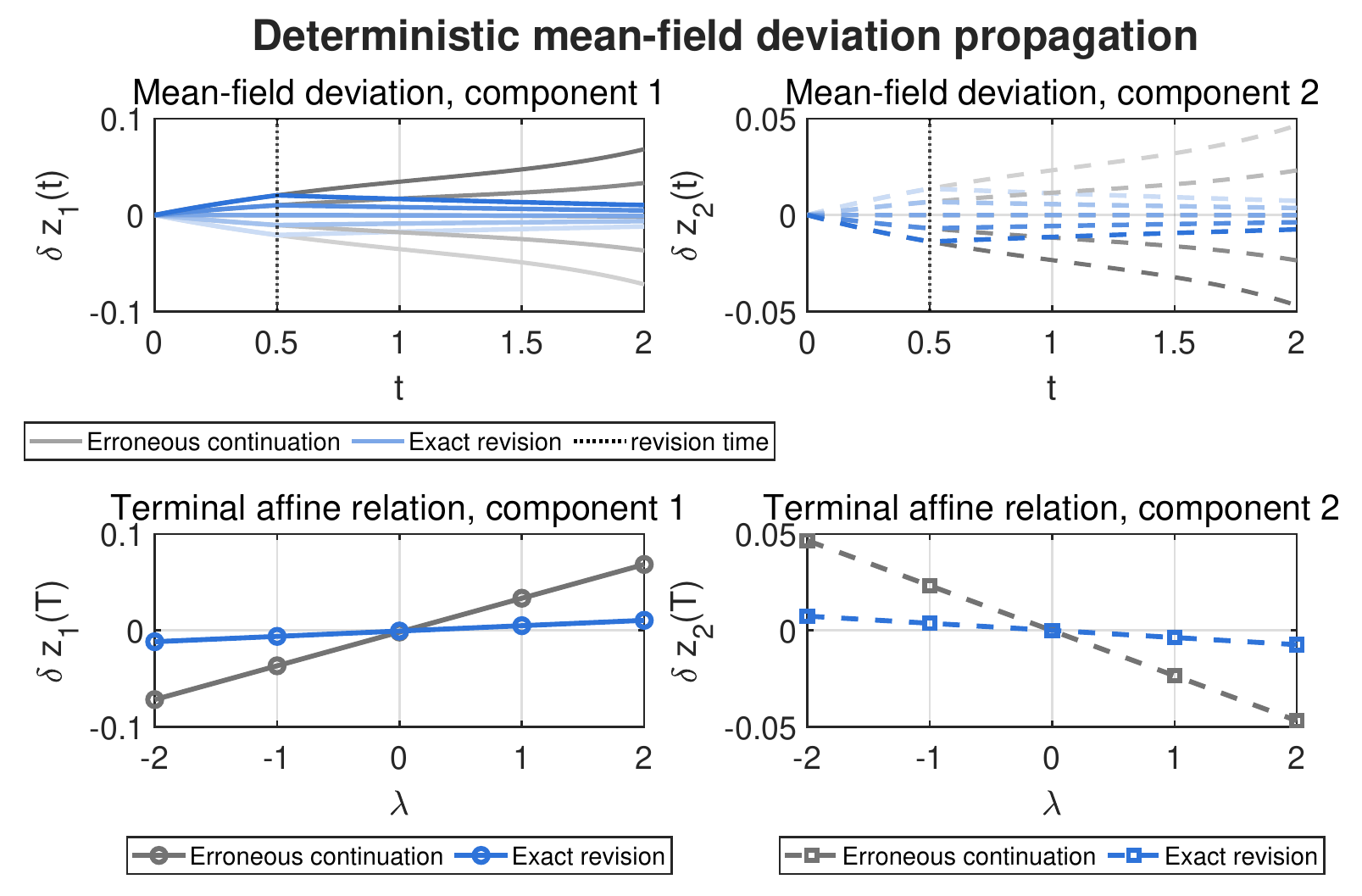}
    \caption{Linear propagation of the actual mean field. The upper row plots the deterministic
    components of \(\delta z_t^A\) and the exact revision
    deviation for \(\bar E=\lambda E_\star\). The lower row shows the
    terminal value as a linear function of \(\lambda\).}
    \label{fig:error_channel}
\end{figure}

This experiment shows the part driven by \(\bar E\) in
Theorem~\ref{thm:sec3_linear_propagation}. Set
\[
\begin{aligned}
\bar E=\lambda E_\star,E_\star=(0.3,-0.2)^\top,\lambda\in\{-2,-1,0,1,2\}.
\end{aligned}
\]
Figure~\ref{fig:error_channel} shows the deterministic mean field deviation
\(\delta z_t^A=z_t^A-z_t^c\) and the deviation after revision.
The terminal panels show the linear relationship in \(\lambda\), matching
\(\delta z_t^A=M_t^z\bar E\).
For the parameter instance considered, exact revision reduces
the displayed terminal deviation.

\subsection{Finite population approximation}

To illustrate Proposition~\ref{prop:sec3_finite_population}, we sample
i.i.d. copies of the generic background triple \((X_0,E,W)\), assign them to
the finite-population coordinates
\((x_0^{j,N},E_j^N,W^{j,N})\), and use the resulting empirical means directly.
For each
\[
N\in\{50,100,200,500,1000,2000\},
\]
we use \(M=20\) Monte Carlo repetitions. For repetition \(m\) and grid time
\(t_k\), define
\[
e_{m,k,N}:=
\left\|x_{t_k}^{(N),A,m}-z_{t_k}^A\right\|^2 .
\]
The reported statistic is
\[
\widehat e_N^{\mathrm{grid}}:=\max_{0\le k\le 400}\frac{1}{M}\sum_{m=1}^M e_{m,k,N}.
\]
This statistic first averages over the Monte Carlo repetitions and then
maximizes over grid time; it is the grid analogue of
\(\sup_{0\le t\le T}\mathbb E|x_t^{(N),A}-z_t^A|^2\) in
Proposition~\ref{prop:sec3_finite_population}.

\begin{figure}[t]
    \centering
    \includegraphics[width=\linewidth]{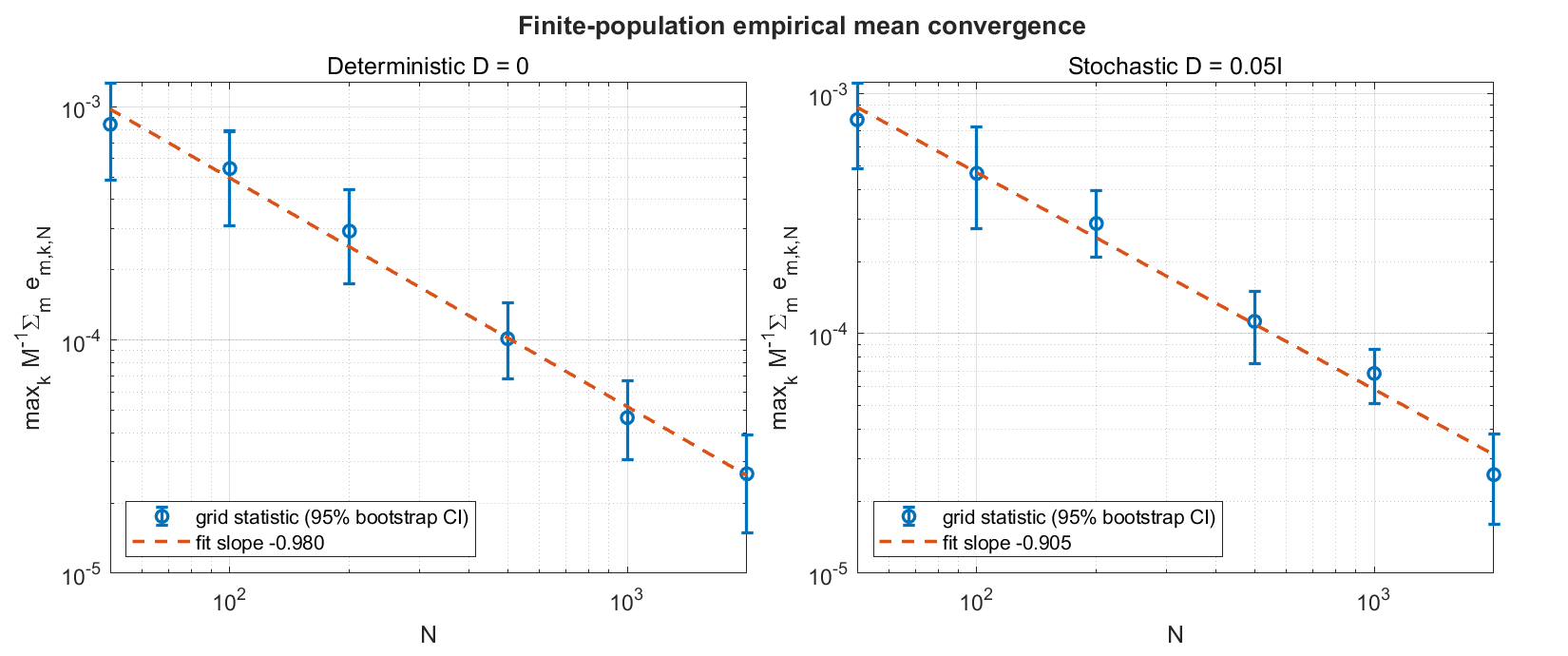}
    \caption{Finite population approximation under the erroneous
    continuation. The plotted quantity is
    \(\max_k M^{-1}\sum_m
    \|x_{t_k}^{(N),A,m}-z_{t_k}^A\|^2\); bars show 95\% bootstrap
    intervals computed from \(500\) resamples.}
    \label{fig:finiteN_convergence}
\end{figure}

Figure~\ref{fig:finiteN_convergence} shows the relation between
\(\widehat e_N^{\mathrm{grid}}\) and \(N\). The fitted slopes and the scaled
quantity \(N\widehat e_N^{\mathrm{grid}}\), reported in
Table~\ref{tab:finiteN_convergence}, are consistent with the \(O(N^{-1})\)
estimate in Proposition~\ref{prop:sec3_finite_population}.

\begin{table}[t]
    \centering
    \scriptsize
    \caption{Finite population convergence summary with \(20\) Monte Carlo repetitions.}
    \label{tab:finiteN_convergence}
    \begin{tabular}{lcc}
        \hline
        Case & fitted slope & range of scaled MSE \\
        \hline
        \(D=0\) & \(-0.980\) & \(0.042\)--\(0.059\) \\
        \(D=0.05I_2\) & \(-0.905\) & \(0.039\)--\(0.068\) \\
        \hline
\end{tabular}
\end{table}

\subsection{Deterministic exact revision and continuation cost}
To quantify the effect of revision, consider \(\mathcal A_i\) with
\[
\begin{aligned}
D&=0,x_0^i=z_0=(0.5,-0.5)^\top,\\
E_i&=(0.3,-0.2)^\top.
\end{aligned}
\]
The general and centered cases take \(\bar E=(0.5,-0.5)^\top\) and
\(\bar E=0\), respectively.
\begin{table}[H]    
    \centering
    \scriptsize
    \caption{Continuation costs for a fixed agent without revision
    (\(J_{i,t_0}^A\)) and with exact revision (\(J_{i,t_0}^r\)) at
    \(t_0=0.5\) on the \(400\)-step grid.}
    \label{tab:main_summary}
    \begin{tabular}{@{}lccc@{}}
        \hline
        Case & $J_{i,t_0}^A$ & $J_{i,t_0}^r$
        & $J_{i,t_0}^A-J_{i,t_0}^r$ \\
        \hline
General instance & 0.330456 & 0.303635 & 0.026821 \\
Centered case & 0.293023 & 0.291881 & 0.001143 \\
\hline

    \end{tabular}
\end{table}
The sign of the general-case gap is specific to the parameter instance.
In the centered case, Corollary~\ref{cor:sec4_centered_performance_gain}
gives a nonnegative continuous-time gap, consistent with the grid comparison.
\subsection{Single-layer signal-based revision}
\label{subsec:numerics_single_layer}
\newcommand{\stocostmcruns}{1000}
\newcommand{\stocostseed}{9400}
\newcommand{\stocostgapnominal}{0.023733}
\newcommand{\stocostcilowernominal}{0.023540}
\newcommand{\stocostciuppernominal}{0.023926}
\newcommand{\stocostnomodification}{0.347475}
\newcommand{\stocostnoupdate}{0.323742}

\rev{The stochastic experiment uses the diagonal single-layer case of
Subsection~\ref{subsec:sec5_single_layer_special}. We consider the synthetic,
exogenous signal family
\[
\bar\zeta_i=\hat\zeta_i=\zeta_i^{(\rho)},
\zeta_i^{(\rho)}
=z_{t_0}^A+\rho\bigl(\bar\varepsilon+\widetilde\varepsilon_i\bigr),
\]
where
\(\widetilde\varepsilon_i\sim\mathcal N(0,\Sigma^\varepsilon)\) and
\(\rho\in\{0,0.25,0.5,1,1.5,2\}\). The signal error is independent of
\((X_0,E,W)\) and of the Brownian increments after \(t_0\).
Thus \(\rho=1\) gives the nominal
signal scale, while \(\rho=0\) supplies the exact current mean-field state. Each
realized signal determines one signal-based feedback law at \(t_0\), denoted by
\((x^{i,\nu,\rho},u^{i,\nu,\rho},z^{\nu,\rho},\bar u^{\nu,\rho})\).}

\rev{At \(\rho=0\), both signal components equal \(z_{t_0}^A\), and the
restricted equations coincide with the oracle equations:}
\[z^{\nu,0}=z^o.\]
\rev{The cost gap is evaluated with common random draws across signal scales
and by trapezoidal quadrature on the \(400\)-step Euler--Maruyama grid.}

\begin{figure}[t]
    
    \centering
    \includegraphics[width=\linewidth]{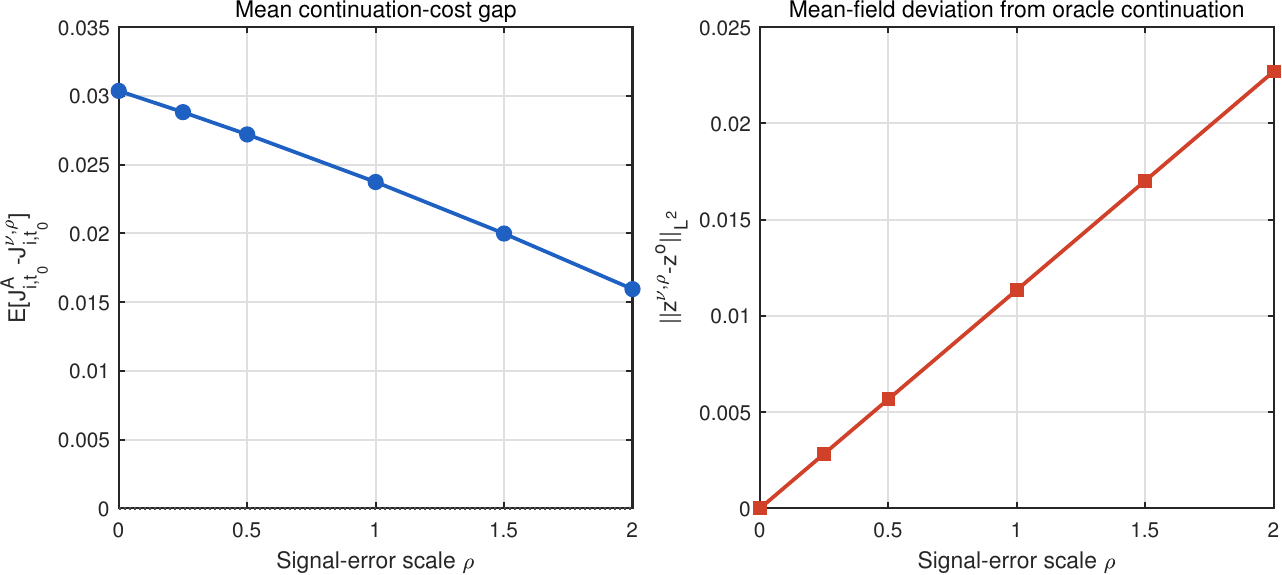}
    \caption{Diagonal single-layer signal-based revision as the
    signal error scale \(\rho\) varies. The left panel reports the mean cost
    gap \(\widehat\Delta J(\rho)\) over \(1000\) paired simulations;
    the right panel reports
    \(\|z^{\nu,\rho}-z^o\|_{L^2(t_0,T)}\), the \(L^2\) deviation
    norm relative to the stochastic oracle continuation initialized at the current state.}
    \label{fig:stochastic_cost_sensitivity}
\end{figure}

At the nominal level \(\rho=1\), the estimated mean continuation costs are
\(\stocostnomodification\) under the erroneous continuation and
\(\stocostnoupdate\) under the signal-based continuation. The paired gap is
\(\stocostgapnominal\). At this parameter value, the estimated mean cost under
the signal-based continuation is lower. As \(\rho\) increases, the plotted
mean-cost gap decreases and the plotted deviation norm increases.
The corresponding centered-state covariance is reported in the
Supplementary Material.

\section{Conclusion}\label{sec:conclusion}

\bluerev{We derived exact linear maps from heterogeneous initial mean-field
observation errors to deviations in subjective predictions, the actual mean
field, individual states, and controls. The propagation operators depend only on the model
data and can therefore be computed independently of the error realizations;
for any fixed agent, the errors enter only through \(E_i\) and \(\bar E\). The
same representation gives covariance, cost, and finite-population relations.

In the deterministic model, recovery of \((E_i,\bar E)\) from a continuously
observed private state trajectory becomes a finite-dimensional observability
problem, with exact recovery if and only if the associated Gramian is
nonsingular. The recovered errors reconstruct \(z_{t_0}^A\) and determine the
correct-information continuation from the state actually reached. For
stochastic dynamics, the two-component signal construction quantifies
deviations from an oracle current-state continuation; the numerical experiments
illustrate these relations. Future work will examine recovery stability under
noisy or partial trajectory observations and derive revision signals from
explicit estimators.}

\clearpage
\section*{Supplementary Material}
\setcounter{section}{0}
\setcounter{theorem}{0}
\setcounter{equation}{0}
\renewcommand{\thesection}{S\arabic{section}}
\renewcommand{\theequation}{S\arabic{section}.\arabic{equation}}
\renewcommand{\thetheorem}{S\arabic{section}.\arabic{theorem}}
\renewcommand{\theHsection}{supp.\arabic{section}}
\renewcommand{\theHequation}{supp.\arabic{section}.\arabic{equation}}
\renewcommand{\theHtheorem}{supp.\arabic{section}.\arabic{theorem}}

This Supplementary Material contains proofs, operator definitions, covariance
and quadratic-cost formulas, and numerical details.

\section{Additional results}

\subsection{Finite-population approximation}

The finite-population approximation is stated in
Proposition~\ref{prop:sec3_finite_population} of the main text. Its proof is
given in Subsection~\ref{supp:proof_sec3_finite_population}.

\subsection{Quadratic cost coefficients}
\label{subsec:supp_cost_coefficients}

This subsection gives the explicit coefficients used in the cost identities
of Sections~3 and 4 of the main text. Define the benchmark residuals
\[
\left\{
\begin{aligned}
&\bar r_{i,t}^I:=\mathbb E[x_t^{i,c}-s],
\bar r_{i,t}^Q:=\mathbb E[x_t^{i,c}-\Gamma z_t^c-\eta],\\
&\bar r_{i,T}^I:=\mathbb E[x_T^{i,c}-\bar s],
\bar r_{i,T}^Q:=\mathbb E[x_T^{i,c}-\bar\Gamma z_T^c-\bar\eta],\\
&\bar u_{i,t}^c:=\mathbb E[u_t^{i,c}].
\end{aligned}
\right.
\]
For the cost distortion over \([0,T]\) in Section~3, the coefficients are
\[
\begin{aligned}
&\ell_i:=\int_0^T\big[X_t^\top Q_I\bar r_{i,t}^I
+U_t^\top R\bar u_{i,t}^c\\
&+(X_t-\Gamma Z_t)^\top Q\bar r_{i,t}^Q\big]\mathrm{d}t
+X_T^\top\bar Q_I\bar r_{i,T}^I\\
&+(X_T-\bar\Gamma Z_T)^\top\bar Q\bar r_{i,T}^Q,
\end{aligned}
\]
and
\[
\begin{aligned}
&\mathcal H:=\int_0^T\big[X_t^\top Q_IX_t+U_t^\top RU_t\\
&+(X_t-\Gamma Z_t)^\top Q(X_t-\Gamma Z_t)\big]\mathrm{d}t\\
&+X_T^\top\bar Q_IX_T
+(X_T-\bar\Gamma Z_T)^\top\bar Q(X_T-\bar\Gamma Z_T).
\end{aligned}
\]
For the continuation-cost comparison in Section~4 and
\(k\in\{A,r\}\), the corresponding coefficients are
\[
\begin{aligned}
&\ell_{i,t_0}^k:=\int_{t_0}^T\big[(\mathsf X_t^k)^\top Q_I\bar r_{i,t}^I
+(\mathsf U_t^k)^\top R\bar u_{i,t}^c\\
&+(\mathsf X_t^k-\Gamma\mathsf Z_t^k)^\top
Q\bar r_{i,t}^Q\big]\mathrm{d}t\\
&+(\mathsf X_T^k)^\top\bar Q_I\bar r_{i,T}^I\\
&+(\mathsf X_T^k-\bar\Gamma\mathsf Z_T^k)^\top
\bar Q\bar r_{i,T}^Q,
\end{aligned}
\]
and
\[
\begin{aligned}
&\mathcal H_{t_0}^k:=\int_{t_0}^T\big[(\mathsf X_t^k)^\top Q_I\mathsf X_t^k
+(\mathsf U_t^k)^\top R\mathsf U_t^k\\
&+(\mathsf X_t^k-\Gamma\mathsf Z_t^k)^\top
Q(\mathsf X_t^k-\Gamma\mathsf Z_t^k)\big]\mathrm{d}t\\
&+(\mathsf X_T^k)^\top\bar Q_I\mathsf X_T^k\\
&+(\mathsf X_T^k-\bar\Gamma\mathsf Z_T^k)^\top
\bar Q(\mathsf X_T^k-\bar\Gamma\mathsf Z_T^k).
\end{aligned}
\]

{
For the revised continuation in Theorem~\ref{thm:sec4_performance_gap}, use
\[
\begin{aligned}
&\mathsf X_t^A:=X_t,\mathsf Z_t^A:=Z_t,\mathsf U_t^A:=U_t,\\
&\mathsf X_t^{r,E}:=\Psi_x(t,t_0)M_{t_0}^{x,1},\\
&\mathsf X_t^{r,\bar E}:=
\Psi_1(t,t_0)M_{t_0}^z-\Psi_x(t,t_0)M_{t_0}^{x,1},\\
&\mathsf X_t^r:=\begin{pmatrix}\mathsf X_t^{r,E}&\mathsf X_t^{r,\bar E}\end{pmatrix},
\mathsf Z_t^r:=\begin{pmatrix}0&\Psi_1(t,t_0)M_{t_0}^z\end{pmatrix},\\
&G_t^{E,r}:=\begin{pmatrix}0&(P_t^0-P_t^1)\Psi_1(t,t_0)M_{t_0}^z\end{pmatrix},\\
&\mathsf U_t^r:=-K_B(P_t^1\mathsf X_t^r+G_t^{E,r}).
\end{aligned}
\]
}

\rev{The deviation maps yield the following general cost expansion and centered
specialization.}

{
\begin{lemma}\label{lem:general_cost_expansion}
Under the standing LQ conditions and
Assumptions~\ref{ass:sec3_probabilistic}--\ref{ass:sec3_subjective_response},
for a fixed private error \(E_i\) and population-average error \(\bar E\),
\[
J_i^A-J_i^c
=\ell_i^\top\theta_i+\frac12\theta_i^\top\mathcal H\theta_i .
\]
The matrix \(\mathcal H\) is symmetric positive semidefinite.
\end{lemma}

\begin{proof}
Apply
\[
\frac12\bigl(\|a+b\|_M^2-\|a\|_M^2\bigr)
=a^\top Mb+\frac12b^\top Mb
\]
to the running and terminal terms, with
\(b=X_t\theta_i\), \(U_t\theta_i\), or
\((X_t-\Gamma Z_t)\theta_i\), and take expectations. This gives
\(\ell_i\) and \(\mathcal H\); positive semidefiniteness follows by applying
the quadratic part to an arbitrary test vector.
\end{proof}

\begin{corollary}\label{cor:sec3_zero_mean_cost}
If \(\bar E=0\), then \(z_t^A=z_t^c\) and
\(\bar u_t^A=\bar u_t^c\). If \(\mathcal H_{11}\) is the upper-left
\(n\times n\) block of \(\mathcal H\), then
\begin{equation}\label{eq:sec3_zero_mean_quadratic_identity}
J_i^A-J_i^c=\frac12E_i^\top\mathcal H_{11}E_i.
\end{equation}
\end{corollary}
}

\subsection{Derivative and sampled recovery}
\label{subsec:supp_recoverability_variants}

We record derivative and sampled variants of the integral observation theorem.
Assume \(D=0\) on \([0,t_0]\), and set
\[
\xi_i:=\begin{pmatrix}E_i\\\bar E\end{pmatrix}.
\]

All expectations in this subsection are under the objective background law:
\[
\bar g_t^A:=\mathbb E[g_t^E]=g_t^c+M_t^g\bar E.
\]

Let \(\Psi_x(t,s):=\Phi_t^x(\Phi_s^x)^{-1}\) be the closed-loop transition
matrix used in Section~4 of the main text.
If a differentiable trajectory is available, define
\[
O_t:=\dot x_t^{i,A}-Ax_t^{i,A}-Bu_t^{i,A}
\]
and
\[
\mathcal R_t^i:=O_t-(C-FK_BP_t^1)z_t^i+FK_Bg_t^i .
\]

\begin{lemma}\label{lem:sec4_residual}
For every \(t\in[0,t_0]\),
\begin{equation}\label{eq:sec4_residual_linear}
\mathcal R_t^i=K_t\xi_i,
\end{equation}
where

\[
\begin{aligned}
&K_t:=\begin{pmatrix}K_t^E&K_t^{\bar E}\end{pmatrix},
K_t^E:=(FK_BP_t^0-C)\Phi_t^1,\\
&K_t^{\bar E}:=(C-FK_BP_t^1)M_t^z-FK_BM_t^g .
\end{aligned}
\]

\end{lemma}

Define the derivative-residual Gramian
\[
\mathcal W_K(t):=\int_0^tK_s^\top K_sds.
\]
Let \(0\le s_0<s_1<\cdots<s_m\le t_0\). For \(k=0,\ldots,m-1\), define
\[
\begin{aligned}
&Y_k^i:=x_{s_{k+1}}^{i,A}-\Psi_x(s_{k+1},s_k)x_{s_k}^{i,A}\\
&-\int_{s_k}^{s_{k+1}}\Psi_x(s_{k+1},t)
\{(C-FK_BP_t^1)z_t^i\\
&-(B+F)K_Bg_t^i\}dt.
\end{aligned}
\]
Then
\begin{equation}\label{eq:sec4_interval_residual}
Y_k^i=\left(\int_{s_k}^{s_{k+1}}\Psi_x(s_{k+1},t)K_tdt\right)\xi_i.
\end{equation}
Set
\[
\mathcal K_s:=
\begin{pmatrix}
\int_{s_0}^{s_1}\Psi_x(s_1,t)K_tdt\\
\vdots\\
\int_{s_{m-1}}^{s_m}\Psi_x(s_m,t)K_tdt
\end{pmatrix},
\mathcal Y_i:=
\begin{pmatrix}
Y_0^i\\
\vdots\\
Y_{m-1}^i
\end{pmatrix}.
\]

\begin{corollary}\label{cor:sec4_integral_sampled_rank}
Assume that
\[
\operatorname{rank}\mathcal K_s=2n .
\]
Then \(\mathcal A_i\) can uniquely recover the private error \(E_i\) and the
population-average error \(\bar E\) from the sampled residuals
\(\{Y_k^i\}_{k=0}^{m-1}\). In this case,
\[
\begin{aligned}
\xi_i
=(\mathcal K_s^\top\mathcal K_s)^{-1}\mathcal K_s^\top\mathcal Y_i .
\end{aligned}
\]
Moreover, \(\mathcal W_K(t_0)\) is nonsingular.
\end{corollary}

\begin{corollary}\label{cor:sec4_analytic_rank}
Suppose that \(K_t\) is analytic in a neighborhood of \(0\). If, for some
integer \(q\ge0\),
\[
\operatorname{rank}
\begin{pmatrix}
K_0\\
\dot K_0\\
\vdots\\
K_0^{(q)}
\end{pmatrix}
=2n,
\]
then \(\mathcal W_K(t)\) is nonsingular for every sufficiently small \(t>0\).
\end{corollary}

\subsection{Signal-error propagation}
\label{subsec:sec5_signal_error_propagation}

{
Let \(z^\star:=z^A_{t_0}\). For
\(\sigma_i=(\bar\zeta_i,\hat\zeta_i)\), define
\[
e_i:=\begin{pmatrix}
\bar\zeta_i-z^\star\\
\hat\zeta_i-z^\star
\end{pmatrix},
e^{\mathsf S}:=\begin{pmatrix}
\bar\zeta-z^\star\\
\hat\zeta-z^\star
\end{pmatrix},
\bar e:=\mathbb E[e^{\mathsf S}].
\]
The stochastic oracle current-state continuation is initialized at
\(z^\star\), and its individual state is coupled with the signal-based state
through the same Brownian increments after \(t_0\). Write
\[
\begin{aligned}
&\Delta z_t:=z_t^\nu-z_t^o,
\Delta x_t^i:=x_t^{i,\nu}-x_t^{i,o},\\
&\Delta u_t^i:=u_t^{i,\nu}-u_t^{i,o}.
\end{aligned}
\]
For the generic background states, set
\[
y_t^\nu:=X_t^\nu-z_t^\nu,y_t^o:=X_t^o-z_t^o.
\]
and define the single centered-state sensitivity matrix
\begin{equation}
\mathsf Y_t^\nu
:=-\Phi_t^x\int_{t_0}^t
(\Phi_s^x)^{-1}BK_B\mathsf C_sds.
\label{eq:supp_sec5_centered_sensitivity}
\end{equation}
Then
\begin{equation}
y_t^\nu-y_t^o
=\mathsf Y_t^\nu(e^{\mathsf S}-\bar e).
\label{eq:supp_sec5_centered_identity}
\end{equation}

\begin{theorem}
\label{thm:two-layer-second-order-distortion}
Under the assumptions of
Theorem~\ref{thm:sec5_signal_error_propagation}, assume
\(e^{\mathsf S}\in L^2\), and let
\[
\Sigma^e:=\operatorname{Cov}(e^{\mathsf S}),
\Sigma_t^\nu:=\operatorname{Cov}(y_t^\nu),
\Sigma_t^o:=\operatorname{Cov}(y_t^o).
\]
For every \(t\in[t_0,T]\),
\begin{equation}
\begin{aligned}
&\Sigma_t^\nu=\Sigma_t^o+\mathsf Y_t^\nu\Sigma^e(\mathsf Y_t^\nu)^\top\\
&+\operatorname{Cov}(y_t^o,e^{\mathsf S}-\bar e)
(\mathsf Y_t^\nu)^\top\\
&+\mathsf Y_t^\nu
\operatorname{Cov}(e^{\mathsf S}-\bar e,y_t^o).
\end{aligned}
\label{eq:supp_sec5_full_covariance}
\end{equation}
If \(e^{\mathsf S}\) is independent of the pair consisting of the oracle
centered state at \(t_0\) and the post-\(t_0\) Brownian increments, then the
cross terms vanish and
\[
\Sigma_t^\nu
=\Sigma_t^o+\mathsf Y_t^\nu\Sigma^e(\mathsf Y_t^\nu)^\top.
\]

For \(\mathcal A_i\), fix the private signal error \(e_i\), the
population-average signal error \(\bar e\), and \(x_{t_0}^{i,A}\), and
condition on these time-\(t_0\) data. The two continuations share their
post-\(t_0\) Brownian increments, so
\(\Delta x^i,\Delta z,\Delta u^i\) are deterministic. Their exact objective
continuation-cost difference is
\begin{equation}
\begin{aligned}
&J_{i,t_0}^\nu-J_{i,t_0}^o=\int_{t_0}^T\Big[
\mathbb E[x_t^{i,o}-s]^\top Q_I\Delta x_t^i\\
&+\mathbb E[u_t^{i,o}]^\top R\Delta u_t^i
+\frac12\|\Delta u_t^i\|_R^2\\
&+\mathbb E[x_t^{i,o}-\Gamma z_t^o-\eta]^\top Q\\
&\cdot(\Delta x_t^i-\Gamma\Delta z_t)
+\frac12\|\Delta x_t^i\|_{Q_I}^2\\
&+\frac12\|\Delta x_t^i-\Gamma\Delta z_t\|_Q^2\Big]dt\\
&+\mathbb E[x_T^{i,o}-\bar s]^\top\bar Q_I\Delta x_T^i\\
&+\mathbb E[x_T^{i,o}-\bar\Gamma z_T^o-\bar\eta]^\top\bar Q\\
&\cdot(\Delta x_T^i-\bar\Gamma\Delta z_T)
+\frac12\|\Delta x_T^i\|_{\bar Q_I}^2\\
&+\frac12\|\Delta x_T^i-\bar\Gamma\Delta z_T\|_{\bar Q}^2.
\end{aligned}
\label{eq:supp_sec5_cost_difference_direct}
\end{equation}
Because the three deviations are linear in \((e_i,\bar e)\), the cost
difference is quadratic in the signal errors.
\end{theorem}
}

\clearpage

\section{Supplementary numerical details and covariance check}\label{subsec:numerics_cov}
\FloatBarrier
{
\subsection{Deterministic trajectory comparison}
The trajectory panels show the correct-information benchmark, the
erroneous continuation, and exact revision from \(z_{t_0}^A\); the marker
indicates \(t_0=0.5\).

\begin{figure}[t]
    \centering
    \includegraphics[width=\linewidth]{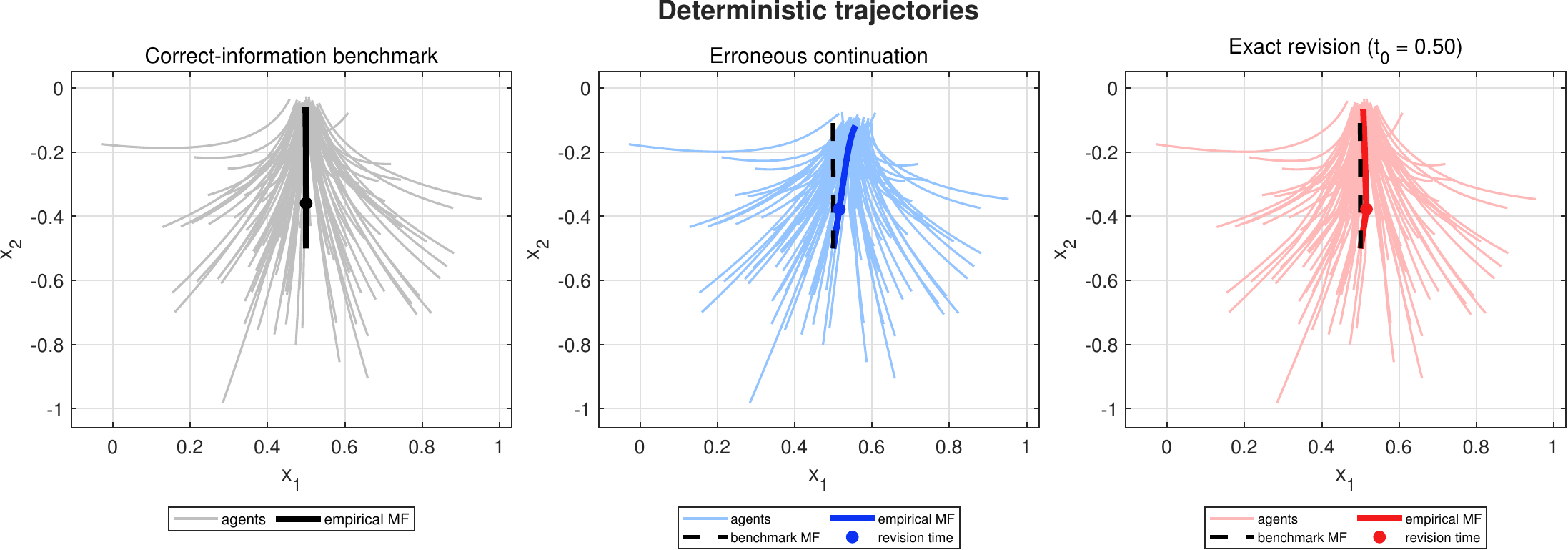}
    \caption{Deterministic trajectory comparison.}
    \label{fig:det_revision_paths}
\end{figure}
}
{
\subsection{Stochastic numerical protocol}
For replication \(m\), let \(L_{i,t_0}^{A,(m)}\) and
\(L_{i,t_0}^{\nu,\rho,(m)}\) be the discretized quadratic losses along the
erroneous and signal-based continuations. For each signal scale, the reported
mean gap is computed from \(1000\) paired simulations as
\[
\widehat\Delta J(\rho)
:=\frac{1}{1000}\sum_{m=1}^{1000}
\left(L_{i,t_0}^{A,(m)}-L_{i,t_0}^{\nu,\rho,(m)}\right).
\]
Within each pair, the two continuations share \(x_0^i\), the private error
\(E_i\), and the post-\(t_0\) Brownian increments; the same signal-error draw
is rescaled across \(\rho\). The
continuation cost in \eqref{eq:sec4_continuation_cost} is evaluated on the
\(400\)-step Euler--Maruyama grid by trapezoidal quadrature.
}
Theorem~\ref{thm:two-layer-second-order-distortion} concerns the centered
generic background states \(y_t^o=X_t^o-z_t^o\) and
\(y_t^\nu=X_t^\nu-z_t^\nu\). For the generic diagonal synthetic signal,
\[
\bar\zeta=\hat\zeta=z_{t_0}^A+\varepsilon,
\]
one has
\(e^{\mathsf S}=(\varepsilon^\top,\varepsilon^\top)^\top\). With
\[
\mathsf Y_t^\varepsilon:=\mathsf Y_t^\nu
\begin{pmatrix}I_n\\I_n\end{pmatrix},
\Sigma^\varepsilon:=\operatorname{Cov}(\varepsilon),
\]
the single-layer covariance formula is
\begin{equation}\label{eq:numerical_diagonal_covariance}
\Sigma_t^\nu
=\Sigma_t^o+\mathsf Y_t^\varepsilon\Sigma^\varepsilon
(\mathsf Y_t^\varepsilon)^\top .
\end{equation}

The experiment uses \(t_0=0.5\), \(D=0.05I_2\),
\(\operatorname{Cov}(X_0)=0.03I_2\),
\(\operatorname{Cov}(E)=0.04I_2\), and
\[
\mathbb E[\varepsilon]=(0.1,-0.1)^\top,
\Sigma^\varepsilon=0.01I_2.
\]
The centered signal error is independent of the centered state under the
stochastic oracle current-state continuation and of all Brownian increments
after \(t_0\). The
nonzero signal mean affects the mean-field path but drops out of the
covariance.

\begin{figure}[t]
    
    \centering
    \includegraphics[width=\linewidth]{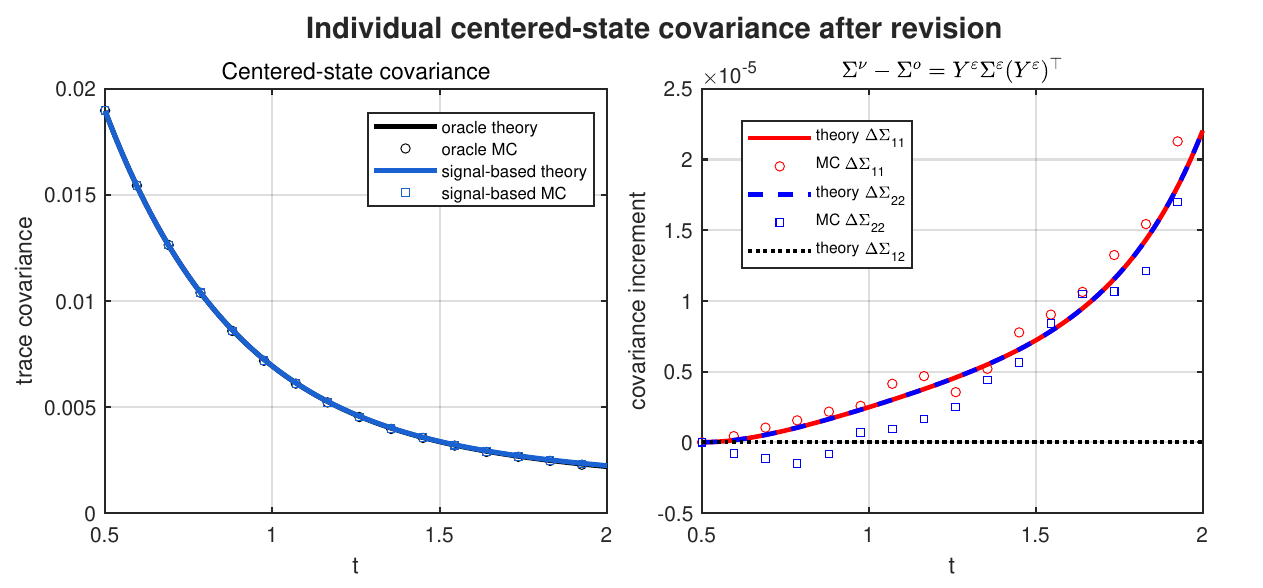}
    \caption{Centered-state covariance. The left
    panel compares theoretical traces with the average of \(80\) batchwise
    cross-sectional estimates, each based on \(N=100\) samples from the
    objective background law; the right panel compares the
    covariance increment
    \(\Sigma_t^\nu-\Sigma_t^o
    =\mathsf Y_t^\varepsilon\Sigma^\varepsilon
    (\mathsf Y_t^\varepsilon)^\top\) with its
    Monte Carlo estimate.}
    \label{fig:cov_prop}
\end{figure}

We use a \(400\)-step Euler grid and \(80\) fixed-seed batches, each containing
\(N=100\) i.i.d. samples from the objective background law. The reported Monte
Carlo covariance is the average of the unbiased cross-sectional covariance
estimates computed in the individual batches.
Figure~\ref{fig:cov_prop} compares these averages with the grid analogues of
\(\Sigma_t^o\) and \eqref{eq:numerical_diagonal_covariance}.  At \(T=2\), the
theory--simulation Frobenius residuals are
\(1.688398\times10^{-5}\) for the stochastic oracle continuation and
\(1.805576\times10^{-5}\) for the signal-based continuation.

\begin{table*}[t]
    
    \centering
    \scriptsize
    \caption{Terminal centered-state covariance:
    theory and the average of \(80\) batchwise estimates, each based on
    \(100\) independent samples from the objective background law.}
    \label{tab:terminal_cov}
    \begin{tabular}{@{}lcccc@{}}
        \hline
        Quantity
        & $\Sigma_T^o$ theory & $\widehat\Sigma_T^o$ MC
        & $\Sigma_T^\nu$ theory & $\widehat\Sigma_T^\nu$ MC \\
        \hline
        $\mathrm{tr}\Sigma$
        & $2.186540\times10^{-3}$ & $2.172281\times10^{-3}$
        & $2.230625\times10^{-3}$ & $2.214047\times10^{-3}$ \\
        $\Sigma_{11}$
        & $1.093270\times10^{-3}$ & $1.089388\times10^{-3}$
        & $1.115312\times10^{-3}$ & $1.112749\times10^{-3}$ \\
        $\Sigma_{22}$
        & $1.093270\times10^{-3}$ & $1.082892\times10^{-3}$
        & $1.115312\times10^{-3}$ & $1.101299\times10^{-3}$ \\
        $\Sigma_{12}$
        & $0$ & $9.008275\times10^{-6}$
        & $0$ & $7.844141\times10^{-6}$ \\
        \hline
    \end{tabular}
\end{table*}


\section{Proofs}\label{supp:proofs}

\subsection{Affine verification of the correct-information benchmark}
\label{supp:benchmark_verification}

We verify individual optimality within
\(\mathfrak U^{\mathrm{fb}}_{[0,T]}(z^c,\bar u^c)\) and mean-field consistency
for the pair \((z^c,\varphi^c)\).

For given deterministic mean field paths \((z,\bar u)\), the representative
LQ problem has state dynamics
\[
dx^i_t=(Ax^i_t+Bu^i_t+Cz_t+F\bar u_t)dt+DdW^i_t .
\]
Admissible deviations are feedback maps in
\(\mathfrak U^{\mathrm{fb}}_{[0,T]}(z,\bar u)\). The stochastic maximum
principle gives
\[
u_t^i=-K_Bp_t^i
\]
and the adjoint equation
\[
\left\{
\begin{aligned}
&dp^i_t=-\big(A^\top p^i_t+(Q_I+Q)x^i_t-Q\Gamma z_t\\
&-Q_Is-Q\eta\big)dt+q^i_tdW^i_t,\\
&p^i_T=(\bar Q_I+\bar Q)x^i_T-\bar Q\bar\Gamma z_T-\bar Q_I\bar s-\bar Q\bar\eta.
\end{aligned}
\right.
\]
Substituting \(p_t^i=P_t^1x_t^i+g_t\) yields
\eqref{eq:sec2_P1}, \eqref{eq:sec2_g}, and the feedback
\eqref{eq:sec2_phi_c}.

Sufficiency follows from the standard LQ verification identity. If
\(\alpha\in\mathfrak U^{\mathrm{fb}}_{[0,T]}(z,\bar u)\) and \(x^\alpha\)
is its closed-loop state, It\^o's formula for the quadratic-affine value
function determined by \(P^1\) and \(g\) gives, with
\(\Delta u_t^\alpha:=\alpha(t,x_t^\alpha)-\varphi(t,x_t^\alpha)\),
\[
\begin{aligned}
&J_i(\alpha;z,\bar u)-J_i(\varphi;z,\bar u)\\
&=\frac12\mathbb E\int_0^T
\|\Delta u_t^\alpha\|_R^2dt\ge0,
\end{aligned}
\]
where \(\varphi(t,x)=-K_B(P_t^1x+g_t)\). Since \(R\succ0\), equality requires
\(\alpha(t,x_t^\alpha)=\varphi(t,x_t^\alpha)\)
\(dt\otimes d\mathbb P\)-almost everywhere and, by closed-loop strong
uniqueness, identifies the minimizing state--control pair.

For the benchmark paths, the averaging identity proved in
Remark~\ref{rem:g_affine_z} gives \(p_t^c=P_t^1z_t^c+g_t^c\). Hence
\[
\begin{aligned}
&\mathbb E[\varphi^c(t,X_t^c)]=-K_B(P_t^1z_t^c+g_t^c)\\
&=-K_Bp_t^c=\bar u_t^c .
\end{aligned}
\]
Taking expectation in the closed-loop state equation then gives the
mean field equation in \eqref{eq:sec2_mf_fbsde}. Since the benchmark
forward--backward system is uniquely solvable under the standing LQ conditions,
the representative mean field is \(z^c\).

\subsection{Proof of Remark~\ref{rem:g_affine_z}}

\begin{proof}

Let \(p^{\mathrm{bg}}\) denote the adjoint corresponding to a generic
background state \(X\). Take expectation in its adjoint equation and set
\[
\bar p_t:=\mathbb E[p_t^{\mathrm{bg}}].
\]
Since \(z_t=\mathbb E[X_t]\), the martingale term has zero expectation and
\(\bar p\) satisfies
\[
\left\{
\begin{aligned}
&\dot{\bar p}_t
=-(A^\top\bar p_t+\mathcal Qz_t-Q_Is-Q\eta),\\
&\bar p_T=\bar{\mathcal Q}z_T-\bar Q_I\bar s-\bar Q\bar\eta .
\end{aligned}
\right.
\]
The same backward equation and terminal condition define the mean field adjoint
in \eqref{eq:sec2_mf_fbsde}. By uniqueness,
\(\bar p_t=p_t=P_t^0z_t+\mathcal G_t\). On the other hand, the affine
representation \(p_t^{\mathrm{bg}}=P_t^1X_t+g_t\) gives
\[
\bar p_t=P_t^1z_t+g_t .
\]

Therefore
\[
g_t=(P_t^0-P_t^1)z_t+\mathcal G_t .
\]
\end{proof}

\subsection{Proof of Theorem~\ref{thm:sec3_linear_propagation}}
\label{supp:proof_sec3_linear_propagation}

\begin{proof}
Subtracting the benchmark mean-field equation from \eqref{eq:sec3_zi} gives
\[
\dot{\delta z}_t^i=(A+C-(B+F)K_BP_t^0)\delta z_t^i,\delta z_0^i=E_i.
\]
Hence \(\delta z_t^i=\Phi_t^1E_i\), and
\eqref{eq:sec3_ubari} gives
\(\delta\bar u_t^i=-K_BP_t^0\Phi_t^1E_i\). Applying
Remark~\ref{rem:g_affine_z} to the subjective and benchmark continuations
yields
\[
\delta g_t^i=(P_t^0-P_t^1)\delta z_t^i=M_t^gE_i .
\]
Applying the same identity to the generic private error \(E\) gives
\(\delta g_t^E=M_t^gE\).

Set \(\delta X_t^A:=X_t^A-X_t^c\). Subtracting the generic benchmark state
equation from \eqref{eq:sec3_background_dynamics} gives
\[
\left\{
\begin{aligned}
&d\delta X_t^A=\big[(A-BK_BP_t^1)\delta X_t^A
-BK_B\delta g_t^E+C\delta z_t^A\\
&+F\delta\bar u_t^A\big]dt,\delta X_0^A=0.
\end{aligned}
\right.
\]
The consistency conditions imply
\[
\begin{aligned}
&\delta z_t^A=\mathbb E[\delta X_t^A],\\
&\delta\bar u_t^A=-K_B\bigl(P_t^1\delta z_t^A
+\mathbb E[\delta g_t^E]\bigr).
\end{aligned}
\]
Since \(\delta g_t^E=M_t^gE\), taking expectation in the generic deviation
equation gives
\[
\left\{
\begin{aligned}
&\dot{\delta z}_t^A=(A+C-(B+F)K_BP_t^1)\delta z_t^A\\
&-(B+F)K_BM_t^g\bar E,\delta z_0^A=0.
\end{aligned}
\right.
\]
Variation of constants therefore gives \(\delta z_t^A=M_t^z\bar E\).
Moreover, \(\mathbb E[\delta g_t^E]=M_t^g\bar E\), and hence
\[
\delta\bar u_t^A=-K_B(P_t^1M_t^z+M_t^g)\bar E .
\]

For \(\mathcal A_i\), the actual deviation satisfies
\[
\left\{
\begin{aligned}
&d\delta x_t^i=\big[(A-BK_BP_t^1)\delta x_t^i
-BK_B\delta g_t^i+C\delta z_t^A\\
&+F\delta\bar u_t^A\big]dt,\delta x_0^i=0.
\end{aligned}
\right.
\]
Let \(r_t^i:=\delta x_t^i-\delta z_t^A\). Then
\[
\dot r_t^i=(A-BK_BP_t^1)r_t^i-BK_BM_t^g(E_i-\bar E),r_0^i=0.
\]
Variation of constants and the definition of \(M_t^{x,1}\) give
\(r_t^i=M_t^{x,1}(E_i-\bar E)\). Adding
\(\delta z_t^A=M_t^z\bar E\) yields
\[
\delta x_t^i=M_t^{x,1}E_i+(M_t^z-M_t^{x,1})\bar E
=M_t^{x,1}E_i+M_t^{x,2}\bar E .
\]
Finally,
\(\delta u_t^i=-K_B(P_t^1\delta x_t^i+\delta g_t^i)\), and the compact
matrix form is exactly \eqref{eq:sec3_error_to_deviation_maps}.
\end{proof}

\subsection{\rev{Proof of Theorem~\ref{thm:sec4_performance_gap}}}

{
For \(k=A\), the deviation identities are those of
\eqref{eq:sec3_error_to_deviation_maps}. For \(k=r\), they follow from
Theorem~\ref{prop:sec4_post_revision_deviation} and
\(g_t^r-g_t^c=G_t^{E,r}\theta_i\). Expanding the quadratic continuation cost
on \([t_0,T]\) around the benchmark gives
\[
J_{i,t_0}^k-J_{i,t_0}^c=(\ell_{i,t_0}^k)^\top\theta_i
+\frac12\theta_i^\top\mathcal H_{t_0}^k\theta_i,k\in\{A,r\}.
\]
Subtracting the two identities yields
\eqref{eq:sec4_exact_performance_gap}.
}

\subsection{Proof of Corollary~\ref{cor:sec3_zero_mean_cost}}

\begin{proof}
The identities \(z^A=z^c\) and \(\bar u^A=\bar u^c\) follow directly from
Theorem~\ref{thm:sec3_linear_propagation} when \(\bar E=0\). Thus the two costs are
evaluated in the same mean field environment, and \(u^{i,c}\) is optimal for
the expected-cost LQ problem with paths \((z^c,\bar u^c)\) within the
admissible feedback class. Applying the completion-of-squares
identity from Subsection~\ref{supp:benchmark_verification} to the erroneous
affine feedback in this common environment gives
\begin{equation}\label{eq:sec3_zero_mean_lq_identity}
\left\{
\begin{aligned}
&J_i^A-J_i^c
=\frac12\mathbb E\Bigl[\int_0^T(\|\delta x_t^i\|_{Q_I+Q}^2+\|\delta u_t^i\|_R^2)dt\\
&+\|\delta x_T^i\|_{\bar Q_I+\bar Q}^2\Bigr].
\end{aligned}
\right.
\end{equation}
For \(\theta_i=(E_i^\top,0)^\top\), this is
\(\frac12E_i^\top\mathcal H_{11}E_i\).
\end{proof}

\subsection{Proof of Corollary~\ref{cor:sec4_integral_sampled_rank}}

\begin{proof}
The actual state equation can be written as
\[
\begin{aligned}
&\dot x_t^{i,A}=(A-BK_BP_t^1)x_t^{i,A}+(C-FK_BP_t^1)z_t^A\\
&-FK_B\bar g_t^A-BK_Bg_t^i.
\end{aligned}
\]
Variation of constants on \([s_k,s_{k+1}]\) gives
\[
\begin{aligned}
&x_{s_{k+1}}^{i,A}=\Psi_x(s_{k+1},s_k)x_{s_k}^{i,A}\\
&+\int_{s_k}^{s_{k+1}}\Psi_x(s_{k+1},t)
\{(C-FK_BP_t^1)z_t^A\\
&-FK_B\bar g_t^A-BK_Bg_t^i\}dt.
\end{aligned}
\]
Subtracting the known reference term with \(z^i,g^i\) gives
\[
\begin{aligned}
&Y_k^i=\int_{s_k}^{s_{k+1}}\Psi_x(s_{k+1},t)
\{(C-FK_BP_t^1)(z_t^A-z_t^i)\\
&-FK_B(\bar g_t^A-g_t^i)\}dt.
\end{aligned}
\]
By Lemma~\ref{lem:sec4_residual}, the integrand in braces is \(K_t\xi_i\).
Thus
\[
Y_k^i
=\left(\int_{s_k}^{s_{k+1}}\Psi_x(s_{k+1},t)K_tdt\right)\xi_i.
\]
Stacking these equations gives
\[
\mathcal Y_i=\mathcal K_s\xi_i .
\]
The stacked matrix has full column rank, so
\(\xi_i=(E_i^\top,\bar E^\top)^\top\) is uniquely determined and
\[
\xi_i=(\mathcal K_s^\top\mathcal K_s)^{-1}\mathcal K_s^\top\mathcal Y_i .
\]

If \(\mathcal W_K(t_0)\) were singular, then some \(\delta\ne0\) would satisfy
\(\int_0^{t_0}|K_t\delta|^2dt=0\). By continuity, \(K_t\delta=0\) on
\([0,t_0]\), and hence
\[
\mathcal K_s\delta = 0.
\]
This contradicts the rank condition. Thus \(\mathcal W_K(t_0)\) is nonsingular.
\end{proof}

\subsection{Proof of Corollary~\ref{cor:sec4_analytic_rank}}

\begin{proof}
Fix \(t>0\) inside the common neighborhood on which \(K\) is analytic. If
\(\mathcal W_K(t)\) were singular, there would exist \(v\ne0\) such that
\[
v^\top\mathcal W_K(t)v=\int_0^t|K_sv|^2\mathrm{d}s=0.
\]
Continuity gives \(K_sv=0\) on \([0,t]\). Analyticity then implies
\(K_0^{(r)}v=0\) for every \(r\ge0\), contradicting the full column rank of
the stack \((K_0^\top,\dot K_0^\top,\ldots,(K_0^{(q)})^\top)^\top\).
Therefore \(\mathcal W_K(t)\) is nonsingular for every sufficiently small
\(t>0\).
\end{proof}

\subsection{Proof of Proposition~\ref{prop:sec3_covariance}}

We use the notation \(\bar E\) and \(\Sigma^E\) introduced in
Assumption~\ref{ass:sec3_probabilistic} of the main text.

\begin{proof}
Set
\[
y_t^A:=X_t^A-z_t^A,y_t^c:=X_t^c-z_t^c .
\]
Subtracting the mean field equations from the generic background state equations
and using Theorem~\ref{thm:sec3_linear_propagation} gives
\[
\left\{
\begin{aligned}
&dy_t^A=[(A-BK_BP_t^1)y_t^A-BK_BM_t^g(E-\bar E)]dt\\
&+DdW_t,y_0^A=X_0-z_0.
\end{aligned}
\right.
\]
and
\[
\left\{
\begin{aligned}
&dy_t^c=(A-BK_BP_t^1)y_t^cdt+DdW_t,\\
&y_0^c=y_0^A.
\end{aligned}
\right.
\]
By variation of constants,
\[
\begin{aligned}
&y_t^A=\Phi_t^x y_0^A
+M_t^{x,1}(E-\bar E)\\
&+\Phi_t^x\int_0^t(\Phi_s^x)^{-1}DdW_s,
\end{aligned}
\]
and
\[
y_t^c
=\Phi_t^x y_0^A
+\Phi_t^x\int_0^t(\Phi_s^x)^{-1}DdW_s.
\]
Hence
\[
\Sigma_t^c = \mathbb E[y_t^c (y_t^c)^\top].
\]
Expanding the product gives
\begin{equation}
\begin{aligned}
&\Sigma_t^c=
\Phi_t^x
\mathbb E[y_0^A (y_0^A)^\top]
(\Phi_t^x)^\top
+\Phi_t^x
\mathbb E[
(\int_0^t
 (\Phi_s^x)^{-1} DdW_s
)\\
&(\int_0^t
 (\Phi_s^x)^{-1} DdW_s
)^\top](\Phi_t^x)^\top,
\end{aligned}
\end{equation}
Using
\[
\mathbb E[y_0^A (y_0^A)^\top] = \Sigma_0,
\]
and the It\^o isometry,
we obtain \eqref{eq:sec3_prop_sigmac}.

Therefore,
\[
\Sigma_t^A = \mathbb E[y_t^A (y_t^A)^\top].
\]
Since
\[
M_t^{x,1}
\mathbb E[(E - \bar E)(E - \bar E)^\top]
(M_t^{x,1})^\top
=
M_t^{x,1} \Sigma^E (M_t^{x,1})^\top,
\]
and the cross terms involving \(E-\bar E\) vanish by the independence of
\(E\) from \(X_0\) and \(W\), together with
\(\mathbb E[E-\bar E]=0\),
we have
\[
\Sigma_t^A
=
\Sigma_t^c
+
M_t^{x,1} \Sigma^E (M_t^{x,1})^\top.
\]
\end{proof}

\subsection{Proof of Proposition~\ref{prop:sec3_finite_population}}
\label{supp:proof_sec3_finite_population}

\begin{proof}
Let \(\Delta_t^{(N)}:=x_t^{(N),A}-z_t^A\). Averaging the finite-population
dynamics under the erroneous continuation and subtracting the mean field
equation for \(z^A\) gives
\[
\left\{
\begin{aligned}
&d\Delta_t^{(N)}=[A+C-(B+F)K_BP_t^1]\Delta_t^{(N)}dt\\
&-(B+F)K_BM_t^g(E^{(N)}-\bar E)dt\\
&+DdW_t^{(N)},\Delta_0^{(N)}=x_0^{(N)}-z_0.
\end{aligned}
\right.
\]
By the variation-of-constants formula,
\begin{equation}
\begin{split}
&x_t^{(N),A}-z_t^A=
\Phi_t^z(x_0^{(N)}-z_0)\\
&+\Phi_t^z\int_0^t(\Phi_s^z)^{-1}[-(B+F)K_BM_s^g(E^{(N)}-\bar E)]ds\\
&+\Phi_t^z\int_0^t(\Phi_s^z)^{-1}DdW_s^{(N)}.
\end{split}
\end{equation}
Using the definition of $M_t^z$, we have
\[
\begin{aligned}
&x_t^{(N),A}-z_t^A=\Phi_t^z(x_0^{(N)}-z_0)+M_t^z(E^{(N)}-\bar E)\\
&+\Phi_t^z\int_0^t(\Phi_s^z)^{-1}DdW_s^{(N)}.
\end{aligned}
\]
Under the i.i.d. hypothesis in Proposition~\ref{prop:sec3_finite_population},
\[
\mathbb E|x_0^{(N)}-z_0|^2=O(N^{-1}),
\mathbb E|E^{(N)}-\bar E|^2=O(N^{-1}).
\]
Define
\[
\widetilde W_t:=\sqrt N W_t^{(N)}
=\frac1{\sqrt N}\sum_{j=1}^NW_t^{j,N}.
\]
Since \(W^{1,N},\ldots,W^{N,N}\) are independent standard Brownian motions
of the same dimension, \(\widetilde W\) is a standard Brownian motion and
\(W^{(N)}=N^{-1/2}\widetilde W\). Thus the It\^o isometry gives
\[
\sup_{0\le t\le T}
\mathbb E\Bigl|\Phi_t^z\int_0^t(\Phi_s^z)^{-1}DdW_s^{(N)}\Bigr|^2
=O(N^{-1}).
\]
The coefficient matrices are bounded on \([0,T]\), and
\eqref{eq:sec3_finite_ms} follows.
\end{proof}

{\subsection{Triangular closure and signal-error propagation}}

{
\begin{proof}
The same affine decoupling identity applies to both average layers. For the
augmented layer, it is
\[
\bar g_t=(P_t^0-P_t^1)\bar z_t+\mathcal G_t.
\]
For the restricted layer associated with a fixed signal \(\sigma\), it is
\begin{equation}
\bar g_t^\sigma=(P_t^0-P_t^1)\bar z_t^\sigma+\mathcal G_t.
\label{eq:supp_sec5_average_decoupling}
\end{equation}
Substitution into
\eqref{eq:sec5_restricted_average_environment} reduces the restricted average
layer to the forward equation
\[
\dot{\bar z}_t^\sigma
=\bigl(A+C-M_BP_t^0\bigr)\bar z_t^\sigma-M_B\mathcal G_t.
\]
Consequently, if \(e=(e^{\bar z},e^{\hat z})\in\mathbb R^{2n}\) denotes an
arbitrary signal error relative to \(z^\star\), then
\begin{equation}
\begin{aligned}
&\bar z_t^\sigma-z_t^o=\Psi_1(t,t_0)e^{\bar z},\\
&\bar g_t^\sigma-g_t^o=(P_t^0-P_t^1)\Psi_1(t,t_0)e^{\bar z}.
\end{aligned}
\label{eq:supp_sec5_average_deviations}
\end{equation}
Thus the average layer is forward unique. Given that layer, the response state
is also a forward linear equation, and variation of constants yields
\begin{equation}
\begin{aligned}
&\delta z_t(e):=z_t^\sigma-z_t^o
=\Phi_t^z(\Phi_{t_0}^z)^{-1}e^{\hat z}\\
&-\Phi_t^z\int_{t_0}^t(\Phi_s^z)^{-1}M_B
(P_s^0-P_s^1)\Psi_1(s,t_0)e^{\bar z}ds,
\end{aligned}
\label{eq:supp_sec5_response_state_direct}
\end{equation}
with predicted mean-control deviation
\begin{equation}
\delta\bar u_t(e)
:=-K_B\Bigl[P_t^1\delta z_t(e)
+(P_t^0-P_t^1)\Psi_1(t,t_0)e^{\bar z}\Bigr].
\label{eq:supp_sec5_response_control_direct}
\end{equation}

For every \(e\), let \(c_t(e)\) be the unique backward solution of
\begin{equation}
\left\{
\begin{aligned}
&\dot c_t(e)=-\bigl[(A^\top-P_t^1BK_B)c_t(e)
+(P_t^1C-Q\Gamma)\delta z_t(e)\\
&+P_t^1F\delta\bar u_t(e)\bigr],c_T(e)=-\bar Q\bar\Gamma\delta z_T(e).
\end{aligned}
\right.
\label{eq:supp_sec5_adjoint_sensitivity}
\end{equation}
The map \(e\mapsto c_t(e)\) is linear, so there is a unique deterministic
matrix \(\mathsf C_t\in\mathbb R^{n\times 2n}\) such that
\[
c_t(e)=\mathsf C_te.
\]
Equations
\eqref{eq:supp_sec5_average_decoupling}--
\eqref{eq:supp_sec5_adjoint_sensitivity} determine the restricted average
layer, response state, and adjoint successively.

For the actual signal-based population, linearity and
\(e^{\mathsf S}\in L^1\) give
\[
\mathbb E[g_t^{\mathsf S}-g_t^o]=\mathsf C_t\bar e.
\]
Taking expectations in the generic state equation and subtracting the oracle
mean equation therefore yields
\[
\frac{d}{dt}\Delta z_t=L_t^1\Delta z_t-M_B\mathsf C_t\bar e,\Delta z_{t_0}=0.
\]
Variation of constants proves
\eqref{eq:sec5_mean_error_direct}. For the individual state, the same-noise
coupling cancels the diffusion terms and gives
\[
\begin{aligned}
&\frac{d}{dt}\Delta x_t^i=(A-BK_BP_t^1)\Delta x_t^i
+(C-FK_BP_t^1)\Delta z_t\\
&-FK_B\mathsf C_t\bar e-BK_B\mathsf C_te_i,\Delta x_{t_0}^i=0.
\end{aligned}
\]
A second variation-of-constants step proves
\eqref{eq:sec5_state_error_direct}, and subtraction of the affine feedbacks
gives \eqref{eq:sec5_control_error_direct}.

For a generic background draw, subtract the mean-field deviation equation
from the state-deviation equation. With
\(q_t:=(X_t^\nu-z_t^\nu)-(X_t^o-z_t^o)\), one obtains
\[
\dot q_t=(A-BK_BP_t^1)q_t-BK_B\mathsf C_t(e^{\mathsf S}-\bar e),q_{t_0}=0.
\]
Variation of constants gives
\eqref{eq:supp_sec5_centered_sensitivity} and
\eqref{eq:supp_sec5_centered_identity}.
\end{proof}
}

\subsection{Proof of Theorem~\ref{thm:two-layer-second-order-distortion}}

{
\begin{proof}
Equation~\eqref{eq:supp_sec5_centered_identity} and bilinearity of covariance
give
\[
\begin{aligned}
&\operatorname{Cov}(y_t^\nu)=\operatorname{Cov}(y_t^o)
+\mathsf Y_t^\nu\operatorname{Cov}(e^{\mathsf S})
(\mathsf Y_t^\nu)^\top\\
&+\operatorname{Cov}(y_t^o,e^{\mathsf S}-\bar e)
(\mathsf Y_t^\nu)^\top\\
&+\mathsf Y_t^\nu
\operatorname{Cov}(e^{\mathsf S}-\bar e,y_t^o),
\end{aligned}
\]
which is \eqref{eq:supp_sec5_full_covariance}. Under the independence
condition in the theorem, \(y_t^o\) is independent of
\(e^{\mathsf S}-\bar e\), so both cross terms vanish.

For the cost identity, the same-noise coupling and common initial state imply
that \(\Delta x^i,\Delta z,\Delta u^i\) are deterministic conditional on the
fixed time-\(t_0\) data. Applying
\[
\frac12\bigl(\|a+b\|_M^2-\|a\|_M^2\bigr)
=a^\top Mb+\frac12\|b\|_M^2
\]
directly to the three running terms and the two terminal terms of the
continuation functional gives
\eqref{eq:supp_sec5_cost_difference_direct}. The propagation formulas show
that every deviation is linear in \((e_i,\bar e)\), which proves the final
quadratic-dependence statement.
\end{proof}
}

\subsection{Proof of Lemma~\ref{lem:sec4_residual}}

\begin{proof}
Using
\[
z_t^A - z_t^i = M_t^z \bar E - \Phi_t^1 E_i
\]
and
\[
\bar g_t^A - g_t^i = M_t^g \bar E - M_t^g E_i,
\]
we obtain
\[
\begin{aligned}
&\mathcal R_t^i=
(F K_B P_t^0-C)\Phi_t^1 E_i\\
&+\bigl((C-FK_BP_t^1)M_t^z-FK_BM_t^g\bigr)\bar E.
\end{aligned}
\]
The coefficient of \(E_i\) follows from
\(M_t^g=(P_t^0-P_t^1)\Phi_t^1\), which gives
\eqref{eq:sec4_residual_linear}.
\end{proof}

\subsection{Mean consistency of the stochastic oracle current-state continuation}

\begin{proof}
Let \(m_t^o:=\mathbb E[X_t^o]\) and
\(v_t^o:=\mathbb E[\varphi^o(t,X_t^o)]\). Since \(\varphi^o\) is affine,
\[
v_t^o-\bar u_t^o=-K_BP_t^1(m_t^o-z_t^o).
\]
Taking expectation in the oracle state equation and subtracting the equation
for \(z^o\) gives
\[
\begin{aligned}
&\frac{d}{dt}(m_t^o-z_t^o)=(A-BK_BP_t^1)(m_t^o-z_t^o),\\
&m_{t_0}^o-z_{t_0}^o=\mathbb E[X_{t_0}^A]-z_{t_0}^o
=z_{t_0}^A-z^\star=0.
\end{aligned}
\]
Uniqueness yields \(m_t^o=z_t^o\), and the affine feedback identity then
gives \(v_t^o=\bar u_t^o\) on \([t_0,T]\).
\end{proof}

\subsection{Verification of the deterministic revised continuation law}

\begin{proof}
For the Section~4 continuation, \(D=0\) and
\(g_t^r=(P_t^0-P_t^1)z_t^r+\mathcal G_t\) by
Remark~\ref{rem:g_affine_z}.
Let \(\bar x_t^r:=\mathbb E[X_t^r]\) and
\(\bar v_t^r:=\mathbb E[\varphi^r(t,X_t^r)]\). Affineness of the feedback
gives \(\bar v_t^r-\bar u_t^r=-K_BP_t^1(\bar x_t^r-z_t^r)\). Averaging the
state equation and subtracting the equation for \(z^r\) gives
\[
\frac{d}{dt}(\bar x_t^r-z_t^r)=(A-BK_BP_t^1)(\bar x_t^r-z_t^r),
\bar x_{t_0}^r-z_{t_0}^r=0 .
\]
Hence \(\bar x_t^r=z_t^r\) and \(\bar v_t^r=\bar u_t^r\) on \([t_0,T]\).
\end{proof}

\end{document}